\documentclass[10pt,a4paper]{article}
\usepackage{indentfirst}
\usepackage{authblk}
\usepackage{amsmath,amssymb,amsfonts,amsthm,graphics,amssymb,amsmath,graphicx,graphics,amsthm,amsfonts,mathrsfs,multicol,multirow, amscd,color}

\usepackage[normalem]{ulem}

\usepackage{amsmath,amsthm}
 \usepackage{amssymb}
 \usepackage{amsfonts}
 \usepackage{mathrsfs}
 \usepackage{amssymb,amsmath,graphicx,graphics,amsthm,amsfonts,mathrsfs,multicol,multirow, amscd,color,amssymb,caption}

\usepackage{graphicx}
\usepackage{makeidx}
\usepackage{xcolor}

\usepackage{enumitem}

\newtheoremstyle{mystyle}
{}
{}
{}
{} 
{\bfseries}
{} 
{.5em}
{}

\theoremstyle{plain}
\newtheorem{theorem}{\bf Theorem}[section]
\newtheorem{proposition}[theorem]{\bf Proposition}

\newtheorem{lemma}[theorem]{\bf Lemma}
\newtheorem{corollary}[theorem]{\bf Corollary}

\usepackage{tikz,chemfig}

\usetikzlibrary{backgrounds}
\usetikzlibrary{shadings,intersections}
\usetikzlibrary{arrows}
\usetikzlibrary{shapes,shapes.geometric,shapes.misc}
\pgfdeclarelayer{edgelayer}
\pgfdeclarelayer{nodelayer}
\usetikzlibrary{automata}
\pgfsetlayers{background,edgelayer,nodelayer,main}
\tikzstyle{none}=[inner sep=0mm]
\tikzstyle{every loop}=[]

\tikzstyle{dottedge}=[dash pattern=on \pgflinewidth off 2pt]
\tikzstyle{dashed}=[dash pattern=on 3pt off 3pt]

\tikzstyle{new style 0}=[fill=black, draw=black, shape=circle]
\tikzstyle{red style 1}=[fill=red, draw=black, shape=circle]
\tikzstyle{blue style 2}=[fill=blue, draw=black, shape=circle]
\tikzstyle{white style 4}=[fill=white, draw=black, shape=circle]
\tikzstyle{bklack style 5}=[fill=black, draw=black, shape=rectangle]
\tikzstyle{red style 3}=[fill=red, draw=black, shape=rectangle]
\tikzstyle{yellow style 7}=[fill=yellow, draw=black, shape=rectangle]
\tikzstyle{new style 8}=[fill={rgb,255: red,0; green,132; blue,0}, draw={rgb,255: red,0; green,131; blue,0}, shape=circle]

\tikzstyle{new edge style 0}=[-]
\tikzstyle{new edge style 1}=[-, draw=red]
\tikzstyle{new edge style 2}=[-, draw=blue]
\tikzstyle{new edge style 3}=[-, draw={rgb,255: red,0; green,156; blue,0}]

\newcommand \equ[2]
{\begin{equation}\label{#1}
		#2
	\end{equation}
}

\newcommand \eqn[2]
{\begin{eqnarray}\label{#1}
		#2
	\end{eqnarray}
}

\def \mato {M^v_1}
\def \matt {M^v_2}

\newcounter{countcase}

\newcounter{countclaim}
\def\inclaim{\addtocounter{countclaim}{1}
	{\vspace{0.2 cm}\noindent {\bf Claim \thecountclaim}: }}

\newcommand{\myproof}{{\noindent {\em Proof}.\quad}\setcounter{countclaim}{0}\setcounter{countcase}{0}}

\newcommand{\proofend}{{\hfill$\Box$}}

\def \setg {{\cal G}}

\def\endproof{\hfill \ \ \hbox{$\sqcup$}\llap{\hbox{$\sqcap$}}}

\newcommand{\al}{\alpha}

\newcommand{\G}{\mathcal{G}}

\def\endproof{\hfill \ \ \hbox{$\sqcup$}\llap{\hbox{$\sqcap$}}}

\newcommand\red[1] {{\color{red} #1}}

\def \setx {{\cal X}}

\begin{document}
	
\title{Characterization of Equimatchable Even-Regular Graphs}

\author[1]{\small Xiao Zhao\thanks{Email:  zhaoxiao05@126.com}}

\author[1]{\small Haojie Zheng\thanks{Email: zhj9536@126.com}}

\author[2]{\small Fengming Dong\thanks{Corresponding author. Emails: fengming.dong@nie.edu.sg and
donggraph@163.com.}}

\author[1]{\small Hengzhe Li\thanks{Email: lihengzhe@htu.edu.cn}}

\author[1]{\small Yingbin Ma\thanks{Email: mayingbincw@htu.cn}}

\affil[1]{\footnotesize
School of Mathematics and Statistics,
Henan Normal University, Xinxiang 453007, P.R. China
}

\affil[2]{\footnotesize National Institute of Education,
	Nanyang Technological University, Singapore}

\date{}

\maketitle

\begin{abstract}
The study of equimatchable graphs traces back to Sumner (1979),	
followed by the foundational work of Lesk, Plummer, and Pulleyblank in 1984. 
Regarding the characterization of equimatchable regular graphs, 
Akbari et al. and Eiben et al. have characterized almost all such graphs,
leaving the class of odd-order  equimatchable $r$-regular graphs with independence number at least $3$
and even $r\geq 6$ unclassified. In the present paper, we give a complete structural characterization of
graphs in this remaining class.
\end{abstract}



\section{Introduction}

All the graphs considered in this paper are connected, finite, simple and undirected.
For any graph $G$, let $V(G)$ and $E(G)$ denote its vertex set and edge set, respectively.
For any $E'\subseteq E(G)$, let
$V(E')$ denote the set of vertices $x\in V(G)$ such that $x$ is an end of some edge in $E'$.
 The \emph{complement} of $G$, denoted by $\overline{G}$,
is the graph with $V(\overline{G})=V(G)$
such that two distinct vertices of
$\overline{G}$ are adjacent if and only if they are not adjacent in $G$.
A {\it matching} of $G$
is a subset $M$ of $E(G)$ such that
each pair of edges in $M$ has no
common end-vertices.
For a  matching $M$ of $G$,
an \emph{$M$-alternating} path in $G$ is a path whose edges are alternately contained
in $M$ and $E(G)\setminus M$.
If neither its origin nor its terminus
of an $M$-alternating path
 is contained in $M$, this path is called an \emph{$M$-augmenting} path.
A matching $M$ in $G$
is called a {\it perfect matching} if every vertex of $G$ is incident with an edge of $M$.
For any  $A\subseteq V(G)$,  let $G[A]$ denote
the subgraph induced by $A$
when $A\ne \emptyset$,
let $G-A$ denote the graph $G[V(G)\setminus A]$
when $A\ne V(G)$,
and let $N_{G}(A)$ denote
the set of vertices each of which is adjacent to at least one vertex of $A$.
For any $v\in V(G)$,
write $N_G(v)$ for $N_G(\{v\})$
and $d_G(v)$ for $|N_G(v)|$.

For any subsets $X$ and $Y$ (not necessarily disjoint) of $V(G)$,
let $E_G(X,Y)$ (or simply $E(X,Y)$) denote the set of edges $xy$ in $G$
with $x\in X$ and $y\in Y$,
and let $e_G(X,Y)=|E_G(X,Y)|$.
We simply write
$E_G(x,Y)$ and $e_G(x,Y)$
for $E_G(\{x\},Y)$
and $e_G(\{x\},Y)$, respectively,
where $x\in V(G)$,
and write $E(X)$ and $e(X)$ for $E_G(X,X)$ and $e_G(X,X)$, respectively. When $Y=V(G)\setminus X$,  $E_G(X,Y)$ is called the {\it edge cut} of $G$ associated with $X$, denoted by $\partial_G(X)$
(or simply
$\partial(X)$).

A graph $G$ is called {\it factor-critical} if $G-v$ has a perfect matching for every $v\in V(G)$.
A graph $G$ is called {\it equimatchable} if all of its maximal matchings have the same size.
For any vertex $v\in V(G)$,
a matching $M$ is said
to {\it isolate} $v$
if $v$ is an isolated vertex in $G-V(M)$.
A matching $M$ isolating a vertex $v$ is called minimal if no subset of $M$ isolates $v$.
The following result due to Eiben and Kotrb\v{c}\'{i}k \cite{Eiben}
characterized the structure of $G-(V(M)\cup \{v\})$  if $v$ is a vertex  isolated
by a minimal matching $M$ in a
$2$-connected factor-critical equimatchable graph $G$.

\begin{theorem}[\cite{Eiben}]
	\label{minimal}
	Let $G$ be a $2$-connected factor-critical equimatchable graph.
	For any vertex $v$ in $G$
	and any minimal matching $M_v$
	isolating $v$,
	$G-(V(M_v)\cup \{v\})$ is isomorphic to $K_{2t}$ or $K_{t,t}$ for some positive integer $t$.
\end{theorem}

The concept of equimatchability was introduced  independently by  Gr\"{u}nbaum \cite{baum}, Lewin \cite{Lewin}, and Meng \cite{Meng} in 1974.
In particular,
Gr\"{u}nbaum \cite{baum} asked for a characterization of
equimatchable graphs.
The first study of this problem 
is the characterization of equimatchable graphs with a perfect matching, that is, 
{\it randomly matchable graphs},
due to Sumner~\cite{Sumner}
in 1979.

\begin{theorem}[\cite{Sumner}]
	\label{sumn}
A nontrivial connected graph $G$ is randomly matchable if and only if $G$ is either $K_{2n}$ or $K_{n,n}$, for $n\geq 1$.
\end{theorem}

In 1984,  a structural characterization of 2-connected equimatchable graphs was provided
by Lesk, Plummer
and Pulleyblank
\cite{Lesk}, as stated below.

\begin{theorem}[\cite{Lesk}]
\label{2-connected}
Any $2$-connected equimatchable graph is either bipartite, or factor-critical or $K_{2t}$ for some $t\geq1$.
\end{theorem}

A subset $S$ of $V(G)$ is called an {\it independent set} of $G$ if
any two vertices in $S$ are not
adjacent in $G$.
The {\it independence number of $G$}, denoted by $\alpha(G)$, is the size of a largest independent set  of $G$.
Due to Eiben and
Kotrb\v{c}\'{i}k \cite{Eiben1},
all odd-order graphs $G$ with   $\alpha(G)\le 2$ are
equimatchable.

\begin{theorem}[\cite{Eiben1}]
	\label{equi}
	Any odd-order graph $G$
	with $\al(G)\le 2$ 
is equimatchable.
\end{theorem}

A consequence follows directly from
Theorems~\ref{2-connected} and~\ref{equi}.

\begin{corollary}\label{equi-co}
Any odd-order $2$-connected graph $G$
with $\al(G)\le 2$  is factor-critical.
\end{corollary}

\myproof
Let $G$ be a
$2$-connected graph
of odd order with $\al(G)\le 2$.
By Theorem~\ref{equi},
$G$ is equimatchable.
Since $|V(G)|$ is odd,
Theorem~\ref{2-connected}
implies that
$G$ is either bipartite or factor-critical.
As $\al(G)\le 2$ and $|V(G)|$ is odd,
if $G$ is bipartite, then
$G$  can only be isomorphic to  $K_{1,2}$ or $K_1$,
contradicting the fact that $G$ is 2-connected. Hence $G$ is factor-critical.
\endproof

Now we focus on the study of
connected and regular graphs which
are equimatchable.
Kawarabayashi,
Plummer and Saito
\cite{Kawarabayashi} showed that the only connected equimatchable 3-regular graphs are $K_4$ and $K_{3,3}$.
Akbari et al. \cite{Akbari}
generalized this result to any connected equimatchable $r$-regular
graph for any odd $r\ge 3$.
The main results 
in \cite{Akbari}
are presented in Theorems~\ref{equi-0-2}
and \ref{equi-0} below.

\begin{theorem}[\cite{Akbari}]\label{equi-0-2}
Any connected equimatchable regular graph is $2$-connected.
\end{theorem}

\begin{theorem}[\cite{Akbari}]
	\label{equi-0}
Let $G$ be 
a connected equimatchable $r$-regular graph.
\vspace{-3mm}
\begin{enumerate}
	[itemsep=-4mm, parsep=0.5cm]
\item If $r$ is odd, then $G$ is either $K_{r+1}$ or 
$K_{r,r}$,  and 

\item if $r$ is even, then
$G$ is $K_{r,r}$ when $|V(G)|$ 
is even, and factor-critical
otherwise, and 

\item in particular, for $r\le 4$, 
$G$ is a graph in the set 
$\{C_3, C_4, C_5, C_7, K_5, K_{4,4}, \overline{C_7},
F_4, \overline{F_4}\}$,
where $F_r$ is the graph in Figure~\ref{figure1-1}.
\end{enumerate}
\end{theorem}

By Theorem~\ref{equi-0} (ii), 
for any connected equimatchable $r$-regular graph $G$, where $r$ is even,
if $|V(G)|$ is even,
then $G\cong  K_{r,r}$;
otherwise, $G$ is factor-critical.
Thus, by Theorems \ref{minimal}, \ref{equi-0-2} and \ref{equi-0} (ii), the following conclusion follows immediately.

\begin{corollary}
	\label{Krr}
Let $G$ be a connected equimatchable $r$-regular graph, where $r$ even. 
If $|V(G)|$ is odd, then 
for any vertex $v$ in $G$
and any minimal matching $M_v$
isolating $v$,
$G-(V(M_v)\cup \{v\})$ is isomorphic to $K_{2t}$ or $K_{t,t}$ for some positive integer $t$.
\end{corollary}

Theorems~\ref{equi} and \ref{equi-0} shows that 
the structure of
the connected equimatchable
$r$-regular graphs $G$ remains unknown for the case that 
$r\ge 6$ is even, $|V(G)|$ is odd
and $\alpha(G)\ge 3$.
In this article, we provide
a structural characterization for such graphs.
A schematic overview of the known results and the open question
on  connected equimatchable
regular graphs
is provided in Figure~\ref{lastpiece}.
\begin{figure}[ht]
	\begin{center}
		\begin{tikzpicture}
			[auto,node distance=2 cm, scale =0.75, transform shape]
			
			\tikzstyle{vertex} = [circle, draw=blue!90]
			
			\tikzstyle{selected vertex} =[rectangle, draw, rounded corners=1.5mm] (rounded) {Rounded};

			\tikzstyle{edge} = [-, very thick]
			
			\node[selected vertex] (v1) {Connected $r$-regular equimatchable graphs $G$};
			
			\node[selected vertex] (v2)
			[below left of=v1]
			{Odd $r$ (Akbari et al. [1])};

			\node[] (v03) [below right of=v1]{};
			
			\node[selected vertex] (v3)
			[right of=v03]
			{\red{Even $r$}};
			
			\node[selected vertex] (v4) [below left of=v2]{$K_{r+1}$};
			
			\node[selected vertex] (v5) [below right of=v2]{$K_{r,r}$};
			
			\node[selected vertex] (v6)
			[below left of=v3]
			{\red{$r\ge 6$}};
			
			\node[] (v07) [below right of=v3]{};
			
			\node[selected vertex] (v7) [right of=v07]{$r\le 4$ (Akbari et al. [1])};

			\node[] (v08) [below left of=v6] {};

\node[selected vertex] (v8)
			[left of=v08]
			{even $|V(G)|$ (Akbari et al. [1])};
			
\node[] (v8-1) [below left of=v8]{};
			
\node[selected vertex] (v8-2)
			[below of=v8]
			{$K_{r,r}$};
			
\path[-] (v8) edge node[left]{} (v8-2);
			
			\node[selected vertex] (v9) [below right of=v6]{\red{Odd $|V(G)|$}};
			
			\node[selected vertex, node distance=2cm] (v10) [below right of=v7]{$C_3, C_4, C_5, C_7, K_5, K_{4,4}, \overline{C_7}, F_4,  \overline{F_4}$};
			
			\node[selected vertex] (v12) [below left of=v9]{$\alpha(G)\le 2$ (\cite{Eiben1})};
			
			\node[] (v013) [below right of=v9]{};
			
			\node[selected vertex] (v13) [right of=v013]{\red{$\alpha(G)\ge 3$ (This paper investigates !!)}};
		
			\path[-] (v9) edge node[left]{} (v12);	
			
			\path[-] (v9) edge node[left]{} (v13);

			\path[-] (v1) edge node[left]{} (v2);
			\path[-] (v1) edge node[below]{} (v3);
			
			\path[-] (v2) edge node[left]{} (v4);
			\path[-] (v2) edge node[below]{} (v5);
			
			\path[-] (v3) edge node[left]{} (v6);
			
			\path[-] (v3) edge node[below]{} (v7);
			
			\path[-] (v6) edge node[left]{} (v8);
			
			\path[-] (v6) edge node[below]{} (v9);
			
			\path[-] (v7) edge node[right]{}
			(v10);
			
		\end{tikzpicture}
	\end{center}

	\caption{A schematic overview of the known results and the open question
	on the study of 
	$r$-regular equimatchable graphs
}
	\label{lastpiece}
\end{figure}

For any even $r\ge 2$,
let $F_r$ denote the graph obtained from the complete bipartite
graph $K_{r,r}$ by deleting all edges
in a matching $M$ of size $\frac{r}{2}$ and adding a
new vertex $u$ and new edges
joining $u$ to all vertices in
$V(M)$, as shown in Figure \ref{figure1-1},
where the dotted lines represent
edges in $M$.

\begin{figure}[ht]
\centering
\includegraphics[width=7 cm] {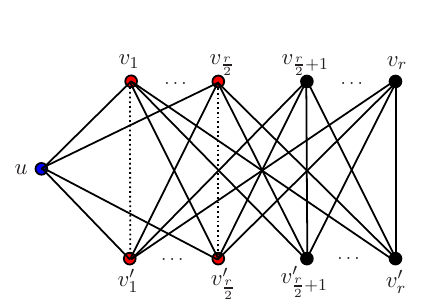}
\caption{Graph $F_r$ for even $r\ge 2$ is 
obtained from  $K_{r,r}$ by deleting all edges
in a matching $M$ in $K_{r,r}$ 
of size $\frac{r}{2}$ and adding a
new vertex $u$ and new edges
joining $u$ to all vertices in
$V(M)$
}
\label{figure1-1}
\end{figure}

In this article, we will prove the following conclusion.

\begin{theorem}\label{main1}
Let $r$ be an even integer with $r\geq 6$. Then every  equimatchable connected
$r$-regular graph $G$ with odd $|V(G)|$ and $\al(G)\geq 3$
is isomorphic to $F_r$.
\end{theorem}

Combining
Theorems~\ref{equi-0}
and \ref{main1},  
we now have a complete 
characterization of equimatchable regular graphs,
as stated in Corollary~\ref{coro0}.

\begin{corollary}\label{coro0}
Let $G$ be a connected equimatchable $r$-regular graph. Then the structure of $G$ is characterized as follows:
\vspace{-3mm}
\begin{enumerate}
		[itemsep=-4mm, parsep=0.5cm]
\item(\cite{Akbari}) if $r$ is odd, then $G\in \{K_{r,r},K_{r+1}\}$;
\item(\cite{Akbari})    if $r\in\{2,4\}$, then  $G\in\{C_3, C_5,C_7,K_5,K_{4,4},\overline{C_7},F_4,\overline{F_4}\}$;
\item(\cite{Akbari}) if $r\geq 6$ is even and $|V(G)|$ is even, then $G\cong K_{r,r}$;
\item(\cite{Eiben1}) if 
$r\ge6$ is even, $|V(G)|$ is odd
and 
$\al(G)\leq 2$,
 then $G$ is always equimatchable; and
\item(This paper) if 
$r\ge6$ is even, $|V(G)|$ is odd
and 
$\al(G)\geq 3$, then $G\cong F_r$.
\end{enumerate}

\end{corollary}

\begin{figure}[ht]
\centering
\includegraphics[width=12 cm] {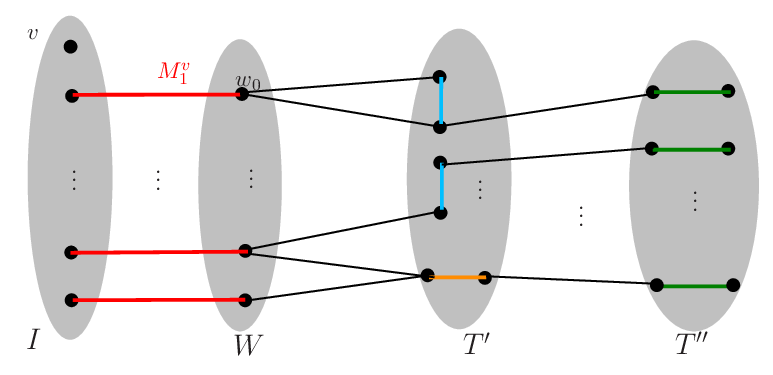}
\caption{$V(G)$ is partitioned into four subsets $I, W, T'$ and $T''$, which are introduced in page~\pageref{nots}.
}
\label{figure2}
\end{figure}

It can be observed 
from Figure~\ref{figure1-1}
that $\{v_1,v_2,\ldots,v_r\}$ is a maximum independent set of $F_r$ of size $r$,
and $F_r-v_1$ has a perfect matching $M^{v_1}$ of size $r$
with $|V(M^{v_1})\cap \{v_1,v_2,\ldots,v_r\}|=r-1$.
Leveraging the properties of $F_r$
mentioned here, 
in order
to prove Theorem \ref{main1},
we decompose the graph $G$
given in Theorem \ref{main1} in the following manner, referring to Figure \ref{figure2}.

Let $r$ be a fixed even integer
with $r\ge 6$
and
$\setg$ be the set of
connected equimatchable $r$-regular graphs $G$
with odd order $|V(G)|$
and  $\alpha(G)\geq 3$.
Let $G=(V,E)\in \setg$.
By Theorem~\ref{equi-0},
$G$ is  factor-critical.
In the following, we introduce
some
notations related to $G$:
$I, v$,
$M^v$, $\mato$, $\matt$,
$M^v_{2,0}$,
$M^v_{2,1}$, $M^v_{2,2}$,
$W$, $T'$, $T''$, $T_w$ for $w\in W$
and $W'$,
which will be applied throughout this article. \label{nots}
\vspace{-3 mm}
\begin{itemize} [itemsep=-1mm]
\item $I$ is an independent set of $G$ with $|I|=\al(G)$
and $v$ is a vertex in $I$;

\item  $M^v$
is a perfect matching of $G-v$, as indicated by the thick lines in Figure \ref{figure2},
which is partitioned into
$\mato$ and $\matt$,
where $\mato$ is the set of edges
$uu'\in M^v$ with $|\{u,u'\}\cap I|=1$, represented by the red thick lines in Figure \ref{figure2}, and $W=V(\mato)\setminus I$, and thus $M^v_2=M^v\backslash M^v_1$;

\item $\matt$ is further partitioned into three subsets
$M^v_{2,0}$, $M^v_{2,1}$ and $M^v_{2,2}$ (represented by green, blue, and orange thick lines, respectively in Figure \ref{figure2}), where, 
 for $i\in \{0,1\}$, 
$M^v_{2,i}$ is the set of edges $uu'\in \matt$ with
$i=\max\{|N_G(u)\cap W|,
|N_G(u')\cap W| \}$. 
Thus, for each $uu'\in M^v_{2,2}$, either 
$|N_G(u)\cap W|\ge 2$ 
or $|N_G(u')\cap W| \}\ge 2$.
\item
$T'=V(M^v_{2,1}\cup M^v_{2,2})$ and
$T''=V(M^v_{2,0})$; and

\item For each $w\in W$,
$T_w$ is the set $\{w\}\cup
\{u,u': uu'\in M^v_{2,1},
w\in N_G(\{u,u'\})\}$,
and $W'$ is the set of $w\in W$
with $|T_w|\ge 3$.
\end{itemize}

	For example, if $G$ is the graph
	$F_r$ shown in Figure~\ref{figure1-1},
	$I$ is the set
	$\{v_1,v_2,\ldots,v_r\}$,
	$v$ is the vertex $v_{r}$
	and $M^v=\{uv'_1\}\cup \{v_iv'_{i+1}: 1\le i\le r-1\}$,
	then $M^v_1=
	\{v_iv'_{i+1}: 1\le i\le r-1\}$,
	$M^v_2=\{uv'_1\}=M^v_{2,2}$,
	$M^v_{2,0}=M^v_{2,1}
	=\emptyset$,
	$W=\{v'_2,\ldots,v'_r\}$,
	$T'=V(M_{2}^v)=V(M_{2,2}^v)=\{u,v'_1\}$ and $T''=\emptyset$.

The structure of our paper is based on the intersection between $M_2^v$ and $W$.
Section~\ref{sec2} establishes foundational structural results for $G\in\mathcal{G}$, relying on Proposition~\ref{move-T} to exclude large odd cliques after deleting a maximum independent set.
Section~\ref{sec3} derives $|W'|=1$ under the assumption of $M_{2,2}^v=\emptyset$. In the Section~\ref{sec4} we apply the conclusion in  Section~\ref{sec3} to show that $M_{2,2}^v\neq\emptyset$.
Section~\ref{sec5} finally verifies $T''=W'=\emptyset$, so $V(G)\setminus\{u\}$ splits into two equal independent sets. By Lemma~\ref{char-Fr}, $G\cong F_r$, completing the proof of Theorem~\ref{main1}.

\section{Properties of graphs in $\setg$
\label{sec2}
}

Recall that $\setg$ is the set of connected equimatchable $r$-regular graphs $G$
	with odd order $|V(G)|$
	and  $\alpha(G)\geq 3$.
	Let $G=(V,E)\in \setg$,
where $r\ge 6$ is even.
In this section, we assume that
$G=(V,E)$ is a graph in $\setg$,
$I$ is an $\alpha$-set of $G$, and $v, M^v, M^v_1, M^v_2,
W, T'$ and $T''$
are as defined on page~\pageref{nots}.

By Theorems~\ref{equi-0-2}
and \ref{equi-0} (ii), 
each graph in $\setg$ is 2-connected
and factor-critical.
Using this fact, we establish the following results.

\begin{proposition}\label{inde}
	For any $G=(V,E)\in \G$
	and maximal matching $M$ of $G$,
	$\alpha(G)
	\le \frac{|V|-1}2$ and
$|M|= \frac{|V|-1}2$.
\end{proposition}

\myproof
Suppose that $\alpha(G)
>\frac{|V|-1}2$.
Then $|I|\ge \frac{|V|+1}2$
for some independent set $I$ of $G$.
Observe that
$$
|\partial(I)|=r\times |I|\ge
\frac{r(|V|+1)}2,
\qquad
|\partial(V\setminus I)|
\le r\times |V\setminus I|
\le \frac{r(|V|-1)}2,
$$
contradicting the fact that
$\partial(I)=\partial(V\setminus I)$.
By Theorem \ref{equi-0} (ii), $G$
is factor-critical,
implying that $|M|=\frac{|V|-1}2$
for some maximal matching.
Since $G$ is also equimatchable,
$|M|= \frac{|V|-1}2$
for each maximal matching $M$. Hence the result holds.
\endproof

\begin{proposition}\label{bipartite}
For any $G=(V,E)\in \G$
and any matching $M$ of $G$,
$\alpha(G-V(M)) \le \frac{|V|-|V(M)|+1}{2}$.
Thus,  if $G-V(M)$ is bipartite
with a bipartition $(A,B)$,
then $-1\le |A|-|B|\le 1$.
\end{proposition}

\myproof Assume that $M$ is a matching
of $G$ such that
$\alpha(G-V(M)) > \frac{|V|-|V(M)|+1}{2}$.
Let $H$ denote the graph $G-V(M)$.
Then, $\alpha(H) =\alpha(G-V(M))> \frac{|V(H)|+1}{2}$,
implying that
$|M'|< \frac{|V(H)|-1}{2}$
for each matching $M'$ of $H$.
It follows that for any maximal matching $M_0$
of $G$ with $M\subseteq M_0$,
$|M_0|< |M|+ \frac{|V(H)|-1}{2}
= \frac{|V|-1}{2}$,
contradicting Proposition~\ref{inde}.

If $G-V(M)$ is  bipartite
	with a bipartition $(A,B)$
	with $|A|\ge |B|$, then
$$
|A|\le \alpha(G-V(M))\le
\frac{|V|-|V(M)|+1}{2}
=
\frac{|A|+|B|+1}2,
$$
implying that $0\le |A|-|B|\le 1$.
The result holds.
\endproof

\begin{lemma}\label{nomatch}
Let $G=(V,E)\in \setg$.
If $X$ is a subset of $V$
such that $T'\cup T'' \subseteq X\subseteq T'\cup T'' \cup W$ and
$|X\cap W|\ge 2$, then $G[X]$ has no perfect matching.
\end{lemma}

\myproof Suppose that
$T'\cup T'' \subseteq X\subseteq T'\cup T'' \cup W$,
$|X\cap W|\ge 2$
and $G[X]$ has a perfect matching $M$.
Note that $|W|=\al(G)-1$ as there is a perfect matching between $W$ and $I-\{v\}$. Then
$|V|- |V(M)|
=|V|-|X|
=|I|+|W|-|X\cap W|
\leq |I|+|W|-2\le 2\al(G)-3$.
It follows that $\al(G)\geq \frac{|V|-|V(M)|+3}{2}$. Note that $V(M)\cap I=\emptyset$,
implying that $I\subseteq V(G-V(M))$.
Therefore $\al(G-V(M))= |I|=\al(G)\geq \frac{|V|-|V(M)|+3}{2}$,
 contradicting  Proposition~\ref{bipartite}.
\proofend

\begin{corollary}
	\label{nomatch-c1}
$W$ is an independent set of $G$
with $|W|=\al(G)-1$.
\end{corollary}

\myproof
Suppose that
$W$ is not an independent set of $G$. Then
	there exists $ww'$ in $E(G)$,
	where $w, w'\in W$.
Let $X=\{w,w'\}\cup V(M^v_{2})$. Obviously,
$T'\cup T'' \subseteq X\subseteq T'\cup T'' \cup W$,
$|X\cap W|=2$, and
$\{ww'\}\cup M^v_{2}$ forms a perfect matching of $G[X]$, which leads to a contradiction
to Lemma~\ref{nomatch}.

By assumption,  $v\in I$
and $M^v$ is a perfect matching of $G-v$.
Since $I$ is an independent set of $G$, we have $|I\cap V(M^v)|
=|I\setminus \{v\}|=|I|-1
=\al(G)-1$,
implying that $|W|=|I\cap V(M^v)|
=\al(G)-1$.
\proofend

\begin{lemma}\label{basic-2}
For any $G=(V,E)\in \setg$,

\vspace{-3 mm}
\begin{enumerate}
	[itemsep=-4mm, parsep=0.5cm]
  \item
  there is no $M_2^v$-augmenting path
  connecting
  two vertices $w_1$ and $w_2$ in $W$;

  \item for any $w\in W$,
  $E_G(T_w, W\setminus \{w\})=\emptyset$;

    \item
   for each $w\in W$,
  $G[T_w]$ is connected with
  $\al(G[T_w])\le 2$ and equimatchable.
  $G[T_w]$ is factor-critical
  whenever it is $2$-connected; and

  \item
  there is at most one vertex
  $u\in T'$ such that
  $|N_G(u)\cap W|\ge 2$, and thus $|M^v_{2,2}|\le 1$.

\end{enumerate}
\end{lemma}

\myproof (i) Suppose that
$P=w_1u_1u_1'\cdots u_{s}u_{s}'w_2$ is an
$\matt$-augmenting path of $G$
connecting two vertices
$w_1$ and $w_2$ in  $W$.
Let $M'$ be the matching
$\{w_1u_1, u_1'u_2,\cdots, u_{s-1}'u_s,
u'_sw_2\}$ and let
$M'':=\matt\setminus E(G[V(P)])$.
Observe that
$M'\cup M''$ is a perfect matching
of $G[T'\cup T''\cup\{w_1,w_2\}]$ and
with $|(T'\cup T''\cup\{w_1,w_2\})\cap W|=2$,
contradicting Lemma~\ref{nomatch}.

(ii). Let $w_1\in W\setminus \{w\}$.
Corollary \ref{nomatch-c1} implies that $w_1w\notin E(G)$.
Suppose that $w_1u\in E(G)$
for some $u\in  T_{w}\setminus \{w\}$.
By the assumptions of $T_w$ and $M_{2,1}^v$, $wu\notin E(G)$
and  $uu'\in M^{v}_{2,1}$ for some vertex
$u'\in N_G(w)$,
implying that $w_1uu'w$ is an $M_2^v$-augmenting path,
contradicting the result in (i).
Thus, the result holds.

(iii) By definition,
$G[T_w]$ is connected.
For any independent set $S$ of $G[T_w]$,  the result in (ii)
implies that $S\cup (W\setminus \{w\})$
is an independent set of $G$.
Thus, $|S|+|W\setminus \{w\}|
=|S|+\al(G)-2\le \al(G)$,
implying that $|S|\le 2$.
Thus, $\al(G[T_w])\le 2$.
By Theorem~\ref{equi},
$G[T_w]$ is equimatchable.
If $G[T_w]$ is $2$-connected, Corollary~\ref{equi-co} implies that $G[T_w]$ is factor-critical
as $|T_w|$ is odd.


(iv) Suppose that
there are $u_1,u_2\in T'$
such that
$|N_G(u_i)\cap W|\ge 2$ for both
$i\in \{1,2\}$.
The result of (i) implies that $u_1u_2\notin M^v_{2,2}$.
Then, $u_1u_1', u_2u'_2\in M^v_{2,2}$
for some vertices $u_1'$ and $u_2'$.
By the result in (i),
$N_G(u_i')\cap W=\emptyset$
for both $i\in \{1,2\}$,
implying that  $W\cup \{u_1',u_2'\}$
is an independent set larger than $I$,
a contradiction.
\proofend

\begin{lemma}\label{basic}
For any $G=(V,E)\in \setg$, the following properties hold:
\vspace{-3 mm}
\begin{enumerate}
	[itemsep=-4mm, parsep=0.5cm]
  \item$N_{G}(u)\cap I\neq \emptyset$
  for each  $u\in V\setminus I$; and

  \item
 $T''\cup(T'\setminus N_{G}(W))$
  	is a clique of $G$,
  implying that
 $T''$ is a clique of $G$
  and $|T''|\leq r$.


\end{enumerate}
\end{lemma}
\myproof
(i). For any $u\in V\setminus I$,
if $N_G(u)\cap I = \emptyset$, then
$I\cup\{u\}$ is an independent set larger than $\alpha=|I|$,
a contradiction.

(ii). Let $V_0:=T''\cup (T'\setminus  N_G(W))$.
Since
$E(T'',W)=\emptyset$,
we have $E_G(V_0,W)=\emptyset$.
Thus,
$$
\al(G)\ge \al(G[W\cup V_0])
=|W|+\al(G[V_0])
=\al(G)-1+\al(G[V_0]),
$$
implying that $\al(G[V_0])\le 1$
and thus $V_0$ is a clique.
Since $T''\subseteq V_0$, 
$G$ is a connected $r$-regular graph, and each vertex in $T''$ has a neighbor in $I$,
 $T''$ must be a clique of $G$
  and $|T''|\leq r$.
\proofend

Now we are going to introduce an
 important result which will be applied repeatedly later.

\begin{proposition}\label{move-T}
	Let $G=(V,E)\in \G$,
	$I'$ be an independent set of $G$
	with $I'\cap I=\emptyset$
	and $|I'|=\al(G)$,
	and $u\in V\setminus (I\cup I')$.
	If $G-(I\cup I'\cup \{u\})$ has
	a perfect matching $M$
	and $|N_G(h)\cap I'|\le 1$
	for each $h\in V\setminus (I\cup I'\cup \{u\})$,
	then
	$G-(I\cup \{u\})$ does not have
	any odd clique $X$
	with the properties that $|X|\ge 3$,
	$|X\cap I'|=1$ and
	$N_G(X)\subseteq \{u\}\cup X\cup I$.
\end{proposition}

\myproof Suppose that
there exists an odd clique $X$ in
$G-(I\cup \{u\})$
such that $|X|\ge 3$,
$|X\cap I'|=1$, and
$N_G(X)\subseteq \{u\}\cup X\cup I$.
We will prove that this assumption  leads to a contradiction.

\begin{figure}[ht]
\centering
\includegraphics[width=9cm] {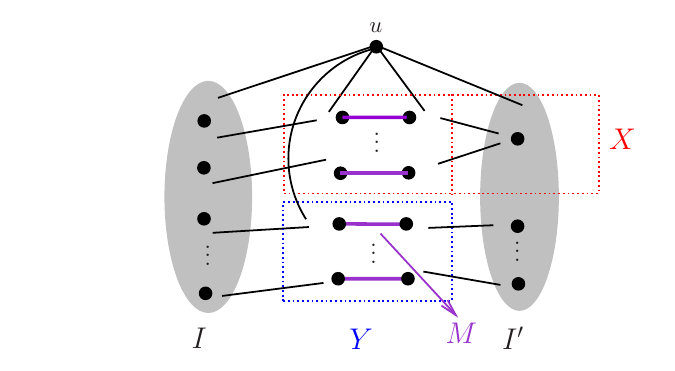}
\caption{$X$ is an odd clique in $G-(I\cup \{u\})$ such that
	$|X|\ge 3$, $|X\cap I'|=1$
and $N_G(X)\subseteq \{u\}\cup X\cup I$ and $Y= V\setminus(I\cup I'\cup X\cup \{u\})$.}
\label{figure3}
\end{figure}

By the assumption of $X$,
$E_G(X,I'\setminus X)=\emptyset$.
Let $Y:= V\setminus(I\cup I'\cup X\cup \{u\})$,
as shown in Figure \ref{figure3}.
Then $E_G(X,Y)=\emptyset$.
Observe that
$|I|=|I'|=\al(G)$ implies that $N_{G}(u)\cap I\ne \emptyset$ and $N_{G}(u)\cap I'\ne \emptyset$.

\inclaim There is no matching
$M_1:=\{uz, x_1z_1,x_2z_2\}$,
where $x_1,x_2\in X$ and $z,z_1,z_2\in I$.

Suppose that such a matching
$M_1$ exists.
Note that
$X\setminus \{x_1,x_2\}$ is an
odd clique,
and for any vertex $x'\in X\setminus \{x_1,x_2\}$,
$G[X\setminus \{x_1,x_2,x'\}]$ has a
perfect matching $M_2$.
If $|X|=3$, then  $M_2=\emptyset$; 
otherwise, a perfect matching in  $G[X\setminus \{x_1,x_2,x'\}]$ can be chosen.
By the given condition,
$G-(I\cup I'\cup \{u\})$ has a perfect matching $M$.
$E_G(X, Y)=\emptyset$ implies that
$M':=M\setminus E(G[X])$
is a perfect matching of
$G[Y]$.
Then,  $M_0:=M_1 \cup M_2
\cup M'$
is a matching of $G$
such  that
$G-V(M_0)$ is bipartite with a bipartition $(I_0$, $I'_0$),
where $I_0=I\setminus \{z,z_1,z_2\}$
and $I'_0=(I'\setminus X)\cup \{x'\}$.
Clearly, 
$|X\cap I'|=1$ and $E_{G}(X,I'\setminus X)=\emptyset$ imply that $|I'_0|=|I'|=|I|$.
Since $|I_0|=\alpha-3=|I'_0|-3$,
it is a contradiction to Proposition~\ref{bipartite}.

\inclaim $|X|=r-1$ and
$|N_G(x)\cap I|
=2-|N_G(u)\cap \{x\}|\le 2$
for each $x\in X$.

Since $X$ is an odd clique and $G$ is $r$-regular, we have $|X|\le r-1$.
Suppose that  $|X|\le r-2$.  Then
$|X|\le r-3$ as $|X|$ is odd.
Thus, $|N_G(x)\cap X|\le r-4$
for each $x\in X$,
implying that $|N_G(x)\cap (V\setminus X)|\ge 4$.
$N_G(X)\subseteq X\cup I\cup \{u\}$
implies that $|N_G(x)\cap I|\ge 3$ for each $x\in X$.
Since $N_G(u)\cap I\ne \emptyset$,
for any $x_1, x_2\in X$,
$G$ has  a matching
$M_1$ with  $M_1\subseteq E_G(I,\{u,x_1,x_2\})$
and $|M_1|=3$,
contradicting Claim 1.
Thus, $|X|=r-1$.

For each $x\in X$, since $X$ is a clique
of size $r-1$ and $N_G(x)\subseteq I\cup X\cup \{u\}$,
we have
$|N_G(x)\cap I|
=2-|N_G(u)\cap \{x\}|\le 2$.
The claim holds.

\inclaim
If $|N_G(u)\cap X|\le r-3$, then
there exists $I_1\subset I$
with $|I_1|=2$ such that
$N_G(x)\cap I\subseteq I_1$
for each $x\in X\cup \{u\}$,
and $G-(X\cup I_1)\cong K_{t,t}$
for some $t\ge \al-1$.

Assume that $|N_G(u)\cap X|\le r-3$.
Since $|X|=r-1$ by Claim 2,
there exist two vertices $x_1,x_2\in X\setminus N_G(u)$.
By Claim 2 again,
$|N_G(x_i)\cap I|=2$
for both $i\in\{1,2\}$.
As $N_{G}(u)\cap I\neq \emptyset$,
	Claim 1 guarantees that
$N_{G}(x_1)\cap I=N_{G}(x_2)\cap I$
and
$N_{G}(u)\cap I\subseteq N_{G}(x_1)\cap I$.
Let $I_1=N_{G}(x_1)\cap I=\{z_1,z_2\}$.
By Claim 1 again, we have $N_{G}(x)\cap I\subseteq I_1$ for each $x\in X$.
Since $X$ is a clique of size $r-1$, $G[X]-\{x_1,x_2,x_3\}$ has
a perfect matching $M_1$ for any $x_3\in X$.
As $N_{G}(x_3)\cap I\subseteq I_1$, we can assume that $x_3z_1\in E(G)$.
Then $M_0:=\{x_2z_2,x_3z_1\}\cup M_1$ is a perfect matching of $G[(X\cup I_1)] -x_1$, which
 is a minimal matching isolating $x_1$.
By Corollary \ref{Krr},
$G_0:=G-(V(M_0)\cup \{x_1\})=G-(X\cup I_1)$
is isomorphic to $K_{2t}$ or $K_{t,t}$
for some $t\ge 1$.
Obviously, $\{u\}\cup(I\setminus I_1)\subseteq V(G_0)$.
As  $N_G(u)\cap I\subseteq I_1$,
$\{u\}\cup(I\setminus I_1)$ is an independent set with
$|\{u\}\cup(I\setminus I_1)|=\al-1\geq 2$ in $G_0$,
implying that $G_0\cong K_{t,t}$,
where $t\ge |\{u\}\cup(I\setminus I_1)|=\al-1$.
Claim 3 holds.

\inclaim
	If there exists a proper subset
	$I_0$ of $I$ such that
either $G-(X\cup I_0)\cong K_{t,t}$
or $G-(X\cup I_0\cup \{u\})\cong K_{t,t}$ for some $t\ge 1$,
then $Y=\emptyset$.

Assume that $V_0:=X\cup I_0$ or   $V_0=X\cup I_0\cup \{u\}$,
where $I_0$ is a proper subset of $I$,
such that $G_0:=G-V_0\cong
K_{t,t}$ for some $t\ge 1$.
Obviously,
$(I\cup I'\cup Y\cup \{u\})\setminus V_0=V(G_0)$.
Suppose that $Y\ne \emptyset$,
i.e.,
$V(G_0)\ne (I\cup I'\cup \{u\})\setminus V_0$.
For any $y\in Y$,
since $|I'|=\al(G)$ and $E_G(X,Y)=\emptyset$,
we have
$E_{G}(I'\setminus X,y)\neq\emptyset$.
Since $G_0\cong K_{t,t}$ and
$I'\setminus X$ is an
independent set of $G_0$,
we have
$I'\setminus X\subseteq N_G(y)$.
This implies that $d_G(y)\ge
| I'\setminus X|=\al-1\ge 2$,
contradicting
 the condition that
 $|N_G(h)\cap I'|\leq 1$
 for each $h\in V\setminus (I\cup I'\cup \{u\})$.
Thus,  Claim 4 holds.

\inclaim $|N_G(u)\cap X|\ge r-2$.

Suppose that $|N_G(u)\cap X|\le r-3$.
By Claim 3, $G-(X\cup I_1)\cong K_{t,t}$,
where $t\ge \al-1$ and
$I_1\subseteq I$ with $|I_1|=2$.
By Claim 4,
$V=I\cup I'\cup X\cup \{u\}$
and $Y=\emptyset$.

Since $V=I\cup I'\cup X\cup \{u\}$
and $G_0=G-(X\cup I_1)$,
$V(G_0)$ can be partitioned into two independent sets:
$(I\setminus I_1)\cup \{u\}$
and $I'\setminus X$,
implying that $G_0\cong K_{\al-1,\al-1}$.
It follows that $I'\setminus X\subseteq N_G(u)$,
and thus
$$
r=d_G(u)\geq|N_G(u)\cap I|
+|I'\setminus X|
\ge 1+\al-1=\al,
$$
implying $\al\le r$.
However, for each $z\in I\setminus I_1$,
since $N_{G}(x)\cap I\subseteq I_1$ for each $x\in X\cup\{u\}$,
we have $E_G(z,X\cup\{u\})=\emptyset$.
Thus,  $Y=\emptyset$ implies
$N_G(z)\subseteq
I'\setminus X$.
Hence,
$
r=d_G(z)\le  |I'\setminus X|=\al-1,
$
implying that $\al\ge r+1$,
a contradiction.
Claim 5 holds.

\inclaim
$I'\not \subseteq N_G(u)$.

Note that $|X\cap I'|=1$.
Assume that $X\cap I'=\{x_0\}$.
Then, $I, I'\setminus \{x_0\}$
and $X$ are pairwise disjoint.

We claim that
$N_{G}(u)\cap I\ne \emptyset$.
Otherwise,  $I\cup \{u\}$ is an
independent set,
contradicting the fact that
$|I|=\al(G)$.

Since $|N_G(u)\cap I|\ge 1$,
and
$|N_G(u)\cap X|\ge r-2$ by Claim 5,
$d_G(u)=r$ implies that
$|N_G(u)\cap (I'\setminus \{x_0\})|\le 1$.
Since $|I'\setminus \{x_0\}|=\alpha-1\ge 2$,
we have $I'\setminus \{x_0\}\not\subseteq N_G(u)$,
implying that $I'\not\subseteq N_G(u)$.
The claim holds.

\inclaim
There exists $z\in I$ such that
$G-(X\cup \{u,z\})\cong K_{t,t}$,
where $t\ge \al-1$.

Let $X_0=N_G(u)\cap X$. By Claim 5,
$|X_0|=|N_G(u)\cap X|\ge r-2\ge 4$.
$|X|=r-1$ and $N_{G}(X)\subseteq I\cup X\cup\{u\}$ imply that $|N_{G}(x)\cap I|=1$ for
each vertex $x$ in $X_0$.
If $N_{G}(x_1)\cap I\neq N_{G}(x_2)\cap I$ for
 any two vertices $x_1, x_2\in X_0$,
 then $|X_0|\geq 4$ implies
 the existence of a
matching $\{uz, x_1z_1, x_2z_2\}$,
where $z,z_1,z_2\in I$,
contradicting Claim 1.
It follows that there are
$x_1,x_2\in X_0$
such that $N_G(x_1)\cap I=N_G(x_2)\cap I=\{z\}$ for some $z\in I$.
Clearly,  $G[(X\cup \{z,u\})]-x_1$
has a  perfect matching $M_1$
with $\{x_2z, x_3u\}\subseteq M_1$,
where $x_3\in X_0\setminus \{x_1,x_2\}$.
Observe that
$M_1$ is a minimal matching isolating $x_1$.
By Corollary \ref{Krr},
$G':=G-(V(M_1)\cup \{x_1\})$
is isomorphic to $K_{2t}$ or $K_{t,t}$
for some $t\ge 1$.
Since $I\setminus \{z\}$ is an independent set in $G'$
of size $\al(G)-1\ge 2$,
we have $G'\cong K_{t,t}$, where
$t\ge \al(G)-1$.

Observe that $G'$ is actually the graph $G-(X\cup \{u,z\})$. The claim holds.

\inclaim $\al=r$.

By Claims 4 and 7, we have
$V=I\cup I'\cup X\cup \{u\}$
and  $G_0:=G-(X\cup \{u,z\})
=G[I\cup I']-\{z,w\}$,
where $w$ is the vertex in $X\cap I'$.
Thus,
$V(G_0)$ can be partitioned into two independent sets:
$I\setminus \{z\}$
and $I'\setminus \{w\}$,
implying that
$G_0\cong K_{\al-1,\al-1}$.
Since $G$ is $r$-regular,
we have $\al-1\leq r$. If $\al-1=r$, then $G_0\cong K_{r,r}$, which implies that $G$ is disconnected.
Thus $\al-1\le r-1$, i.e., $\al\le r$.

By Claim 6, there exists $w'$ in
$I'\setminus N_G(u)$.
Since $N_G(w')\subseteq
V\setminus (I'\cup X\cup\{u\})=I$,
 we have
$r=d_G(w')\le  |I|$,
implying that $\al\ge r$.
Hence $\al=r$.
Claim 8 holds.

Now we are going to complete the proof.
By Claim 2, we have
 $|X|=r-1$ and $|N_G(x')\cap I|=2$
 for each $x'\in X\setminus N_G(u)$.
By Claim 5, $|N_G(u)\cap X|\ge r-2$.  Thus, $|X\setminus N_G(u)|\leq 1$.
It follows that $|N_{G}(X\setminus N_G(u))\cap I|\leq 2$.
Since $G$ is $r$-regular and $|N_G(u)\cap X|\ge r-2$, $|N_{G}(u)\cap I|\le 2$.
Thus,
$|N_{G}(u)\cap I|+
|N_G(X\setminus N_G(u))\cap I|\le 4
<r=\al$,
where the last equality follows from Claim 8,
implying that
$ I\setminus N_G(X\cup \{u\})\neq \emptyset$.
Claims 4 and 7 imply that $Y=\emptyset$ and $V=I\cup I'\cup X\cup \{u\}$.
It follows that $N_G(z')\subseteq I'\setminus X$
for each $z'\in I\setminus
N_G(X\cup \{u\})$. Hence we have
\equ{eq22}
{
r=d_G(z')\le  |I'\setminus X|=\al-1,
}
implying that $\al\ge r+1$,   contradicting Claim 8.
Hence the result follows.
\proofend

We end this section by applying Proposition~\ref{move-T} to prove that
$T'\ne \emptyset$, i.e.,
$N_G(W)\cap V(\matt)\ne \emptyset$.

\begin{lemma}\label{T'non}
For any $G=(V,E)\in \setg$, we have $T'\ne \emptyset$.
\end{lemma}

\myproof
 Suppose that $T'=\emptyset$.
 We claim that
$T''\ne \emptyset$.
Otherwise,
$\al(G)=\frac{|V|+1}{2}$,
contradicting Proposition~\ref{inde}.
Clearly, $|T''|$ is even.
By Lemma~\ref{basic} (ii), $T''$ is a clique.
Assume that $T''\cong K_q$
($q\leq r$).
By the definition of $T''$, we know that $E(W, T'')=\emptyset$.
Therefore, $N_G(W)\subseteq I$
and $N_G(T'')\subseteq I$. Since
$\partial(I)=\partial (W\cup T'')$, it follows that
\equ{eq21}
{
\al(G)\times r =|\partial(I)|
=|\partial(W)|+|\partial(T'')|
=
(\al(G)-1)\times r
+q(r-q+1).
}
Solving  equation (\ref{eq21})
	yields that $q=r$ or $q=1$.
Since $q=|T''|$ is even, we have
$q=r\ge 6$.

Now assume that $u$ and $u'$ are any two vertices in $T''$. Let
$I':=W\cup \{u'\}$ and
$X:=T''\setminus \{u\}$.
Observe that $I'$ is an independent set
with $I'\cap I=\emptyset$ and $|I'|=|I|$ and $|X\cap I'|=1$,
and $G-(I\cup I'\cup \{u\})=G[T''\setminus \{u,u'\}]$
is a complete graph with even order and thus has a perfect matching.
By the definition of $T''$,
we have $N_{G}(h)\cap I'=\{u'\}$ for each $h\in V\setminus(I\cup I'\cup\{u\})= T''\setminus \{u,u'\}$.
But, $X$ is an odd clique of $G-(I\cup \{u\})$ with $|X|=q-1\ge 5$
and $N_G(X)\subseteq I\cup X\cup \{u\}
$, contradicting Proposition~\ref{move-T}.
\proofend

\section{Results under the 
	assumption $M_{2,2}^v=\emptyset$
\label{sec3}
}

In this section,
we establish some conclusions
under the assumption that  $M_{2,2}^v=\emptyset$,
and in the next section,
we will apply these conclusions
to prove that this assumption
leads to a contradiction, and thus
we prove the fact that $M_{2,2}^v\ne \emptyset$.

Suppose that $M_{2,2}^v=\emptyset$ throughout  this section.
It follows that $|N_G(u)\cap W|\le 1$
for each $u\in T'$.
We will show that $|W'|\le 1$,
i.e., there is at most one vertex
$w$ in $W$ with $|T_w|\ge 3$.

Recall that
$M^v_{2,1}$ is the set of edges $u_1u_2\in M^v_2$ such that
$E_G(\{u_1,u_2\}, W)\ne \emptyset$
and $|N_G(u_i)\cap W|\le 1$
for both $i\in\{1,2\}$,
and for any $w\in W$,
$T_w=\{w\}\cup\{u_1,u_2:
u_1u_2\in M^v_{2,1},
w\in N_G(u_1)\cup N_G(u_2)\}$.
Clearly,
$G[T_w]$ is connected for each $w\in W$.
By
Lemma~\ref{basic-2} (ii),
$T_{w_1}\cap T_{w_2}=\emptyset$
for each pair of vertices $w_1,w_2\in W$.
Since $M^v_{2,2}=\emptyset$
by assumption,
$V\setminus (I\cup T'')$ can be partitioned into $T_w$'s for all $w\in W$.

\subsection{Relation between $T_{w_1}$ and $T_{w_2}$}

In this subsection, we mainly
study the relation between
$T_{w}$ and $T''$, or between
$T_{w}$ and $T_{w'}$,
where $w,w'\in W$.

\begin{lemma}\label{basic-3}
Let $G=(V,E)\in \setg$ with
$M_{2,2}^v=\emptyset$.
Then
 the following properties hold:
 \vspace{-3 mm}
\begin{enumerate}
	[itemsep=-4mm, parsep=0.5cm]
  \item for any $w_1, w_2\in W$,
  $G[T_{w_1}\cup T_{w_2}\cup T'']$
  does not have any perfect matching,
  and in particular,
  $G[T_{w_1}\cup T_{w_2}]$ does not have any perfect matching;

  \item for any $w_1, w_2\in W$,
   if both $T_{w_1}$ and $T_{w_2}$ are cliques, then  $E_G(T_{w_1},T_{w_2})=\emptyset$;

  \item  for any $w\in W$,
  	$\al(G[T_{w}\cup T''])\leq 2$;

  	\item  for $X=T_w\cup T''$ or $X=T_w$, where $w\in W$,
  	if  $G[X]$ is connected
  	but not $2$-connected,
  	then
 there exists a cut-vertex $u$ in $G[X]$ such that
 (a) $G[X]-u$ has exactly two components whose vertex sets are
 $O_1$ and $O_2$, 
 (b) both $O_1\cup \{u\}$ and
  	$O_2$ are cliques of $G$,  and  (c) $G[X]-x$ has a perfect matching for any $x\in O_1\cup O_2$.

\end{enumerate}
\end{lemma}

\myproof  By Lemma~\ref{basic-2} (ii),  $E_G(T_{w},W\setminus\{w\})=\emptyset$ for any vertex $w$ of $W$.

(i). Suppose that $M'$ is a perfect matching of $G[T_{w_1}\cup T_{w_2}\cup T'']$.
Let $M'':=M_2^{v}\setminus E(G[T_{w_1}\cup T_{w_2}\cup T''])$.
Then $M'\cup M''$ is a perfect matching of $G[X]$, where
$X:=V(M'\cup M'')$ is a subset of $T'\cup T''\cup W$ with
$|X\cap W|=2$ and $T'\cup T''\subseteq X$,
 contradicting Lemma~\ref{nomatch}.

Now suppose that $M_1$ is a perfect matching of $G[T_{w_1}\cup T_{w_2}]$.
Let $M_1':=M_1\cup M^{v}_{2,0}$.
Then $M_1'$ is a perfect matching of $G[T_{w_1}\cup T_{w_2}\cup T'']$, a contradiction.
Hence (i)  holds.

(ii). Since $T_{w_1}$ and $T_{w_2}$ are both odd cliques, if $E_G(T_{w_1},T_{w_2})\neq\emptyset$,  then  $G[T_{w_1}\cup T_{w_2}]$ has  a perfect matching, contradicting (i).

(iii).
The definition of $T''$ implies that $E(T'', W\setminus\{w\})=\emptyset$. Lemma \ref{basic-2} (ii) implies that $E(T_w, W\setminus\{w\})=\emptyset$.
It follows
that $E(T_w\cup T'', W\setminus\{w\})=\emptyset$.
If $\alpha(G[T_w\cup T''])\ge 3$,
then
\equ{eq24}
{
\al(G)\ge
\al(G[T_w\cup T''])+
\al(G[W\setminus \{w\}])
\ge 3+\al(G) -2=\al(G)+1>\al(G),
}
a contradiction.

(iv). Assume that $G[X]$ is not $2$-connected.
We claim that
$G[X]-u$ has exactly two components
for any cut-vertex $u$ of $G[X]$.
Otherwise,  $\al(G[X])\ge 3$, a contradiction to the conclusion
of (iii) .
Let
$O_1$ and $O_2$ be the vertex sets
of the two components of
$G[X]-u$.
Since $\al(G[X])=2$ and
$u$ is a cut-vertex of $G[X]$,
it can be verified easily that
$O_i$ is a clique
for both $i\in\{1,2\}$.
$G[O_i\cup \{u\}]$ is also a clique
for some $i\in \{1,2\}$.
Otherwise, we can find an independent set of $X$ of size at least 3, contradicting (iii).
Without loss of generality, we assume that $O_1\cup \{u\}$ is a clique.
If both $|O_1|$ and $|O_2|$ are odd or both $|O_1|$ and $|O_2|$ are even and $|N_{G}(u)\cap O_2|\geq 2$, then since both $O_1\cup \{u\}$ and $O_2$ are cliques and $|N_{G}(u)\cap O_2|\geq 1$, $G[X]-x$ has a perfect matching for any $x\in O_1\cup O_2$.
If both $|O_1|$ and $|O_2|$ are even and $|N_{G}(u)\cap O_2|=1$, the only vertex in
$N_{G}(u)\cap O_2$ is a cut-vertex of $G[X]$ which
makes the results hold. (iv) holds.
\proofend

\begin{lemma}\label{basic-4}
Let $G=(V,E)\in \setg$, $\al(G[T_{w}])=2$ holds
  	for at most  one vertex $w\in W$.
  \end{lemma}

\myproof Suppose there are two vertices $w_1$ and $w_2$ in $W$ such that $\al(T_{w_i})=2$ for  both $i\in\{1,2\}$.
We first prove
the following claims.

\inclaim $\al(G[T_{w_1}\cup T_{w_2}])\le 3$.

Suppose Claim 1 fails.
Lemma \ref{basic-2} (ii) implies that $E(T_{w_1}\cup T_{w_2}, W\setminus\{w_1,w_2\})=\emptyset$.
Then
$\al(G[T_{w_1}\cup T_{w_2}])\ge 4$ and
thus,
$$\al(G)\ge\al(G[T_{w_1}\cup T_{w_2}])
+\al(G[W\setminus \{w_1,w_2\}])
= 4+\al(G) -3>\al(G),$$
a contradiction.

\inclaim $E(T_{w_1}, T_{w_2})
\ne \emptyset$.

Suppose that
$E(T_{w_1}, T_{w_2})
= \emptyset$. Then,
$
\al(G[T_{w_1}\cup T_{w_2}])
=\al(G[T_{w_1}])+ \al(G[T_{w_2}])
= 4,
$
a contradiction to Claim 1.
Thus, Claim 2 holds.

\inclaim For any $x_1x_2\in E(T_{w_1},T_{w_2})$,
where $x_i\in T_{w_i}$,
$G[T_{w_i}]-x_i$ has no perfect matchings for some
$i\in \{1,2\}$.

If $G[T_{w_i}]-x_i$ has perfect matchings for both $i\in\{1,2\}$,
then $x_1x_2\in E(T_{w_1},T_{w_2})$
implies that
$G[T_{w_1}\cup T_{w_2}]$
has perfect matchings,
contradicting Lemma~\ref{basic-3} (i).
Thus, $G[T_{w_i}]-x_i$ has no perfect matchings for some
$i\in \{1,2\}$.

\inclaim
$G[T_{w_i}]$ is not 2-connected
for some $i\in \{1,2\}$.

Suppose that both
$G[T_{w_1}]$ and $G[T_{w_2}]$ are 2-connected.
Corollary \ref{equi-co} and $\al(G[T_{w_i}])=2$ imply that $G[T_{w_i}]$ is factor-critical
for both $i\in\{1,2\}$.
By Claim 2,  there exists
$x_1x_2\in E(T_{w_1}, T_{w_2})$.
By Claim 3,
$G[T_{w_i}]-x_i$ has no perfect matchings for some $i\in \{1,2\}$,
implying that
$G[T_{w_i}]$
is not factor-critical,
a contradiction.

\inclaim
$G[T_{w_i}]$ is not 2-connected
for both $i\in\{1,2\}$.

By Claim 4, either
$G[T_{w_1}]$ or $G[T_{w_2}]$
is not 2-connected.
Now suppose that Claim 5 fails,
and $G[T_{w_1}]$ is not 2-connected but
$G[T_{w_2}]$ is 2-connected.
Corollary \ref{equi-co} and $\al(G[T_{w_2}])=2$ imply that $G[T_{w_2}]$ is factor-critical.

Applying Lemma \ref{basic-3} (iv)
to $G[T_{w_1}]$ yields that
 there exists a cut-vertex $u$ in $G[T_{w_1}]$ such that  $G[T_{w_1}]-u$ has exactly two
components whose vertex sets are
$O_1$ and $O_2$
with the properties that both $O_1\cup\{u\}$ and
  	$O_2$ are cliques of $G$ and $G[T_{w_1}]-x$ has perfect matching for each $x\in O_1\cup O_2$.

  	By Claim 2, $E(T_{w_1}, T_{w_2})\ne \emptyset$.
  	If there exists an edge $x_1x_2
  	\in E(O_1\cup O_2, T_{w_2})$,
  	where $x_1\in O_1\cup O_2$
  	and $x_2\in T_{w_2}$,
  	then Claim 3 implies that
  	$G[T_{w_2}]-x_2$ has no perfect
  	matching,
  	contradicting the fact that $G[T_{w_2}]$ is factor-critical.
  	It follows that
  	$E(T_{w_1}, T_{w_2})
  	=E(u, T_{w_2})$,
  	implying that
  	$\al(G[T_{w_1}\cup T_{w_2}])
  	\ge \al(G[(T_{w_1}\setminus \{u\})\cup T_{w_2}])
  	=\al(G[T_{w_1}\setminus \{u\}])
  	+\al(G[T_{w_2}])
  	=4$,
  	contradicting Claim 1.
  	Hence Claim 5 holds.
  	
 Now we are going to complete the proof. By Claim 5, for $i\in\{1,2\}$,
 $G[T_{w_i}]$ is not $2$-connected.
By  Lemma \ref{basic-3} (iv),
for any $i\in\{1,2\}$,
there exists a cut-vertex
$u_i$ in $G[T_{w_i}]$ such that $G[T_{w_i}]-u_i$ has exactly two  components
whose vertex sets are
$O_{i,1}$ and $O_{i,2}$
with the properties that both $O_{i,1}\cup\{u_i\}$ and
$O_{i,2}$ are cliques of $G$ and $G[T_{w_i}]-x$ has perfect matchings for any
$x\in O_{i,1}\cup O_{i,2}$.
Thus, Claim 3 implies that
$E(O_{1,1}\cup O_{1,2},
O_{2,1}\cup O_{2,2})=\emptyset$,
i.e.,
$E(T_{w_1}, T_{w_2})
\subseteq E(u_1, T_{w_2})
\cup E(T_{w_1}, u_2)$.
It follows that
$$
\al(G[T_{w_1}\cup T_{w_2}])
\ge \al(G[T_{w_1}\cup T_{w_2}]-\{u_1,u_2\})
=\al(G[T_{w_1}\setminus\{u_1\}])
+\al(G[T_{w_2}\setminus\{u_2\}])
=2+2=4,
$$
 contradicting Claim 1.
Hence the lemma holds.
\proofend

\subsection{$|T_w|\ge 3$ 	for only one vertex $w\in W$}

Recall that $W'$ denotes the set of vertices $w\in W$ for which $|T_w|\ge 3$.
By Lemma~\ref{T'non}
and the assumption
that  $M_{2,2}^v=\emptyset$
in this section, we have
$W'\ne \emptyset$.

Choose a vertex $w_0$ in $W'$
satisfying the following conditions:
 \vspace{-3 mm}
\begin{enumerate}
	[itemsep=-4mm, parsep=0.5cm]
\item $T_{w_0}$ is not a clique
if there exists some $w\in W'$
such that
$T_w$ is not a clique;
and
\item $|E_G(T_{w_0},T'')|\ge |E_G(T_{w},T'')|$
for each $w\in W'$
if $T_w$ is a clique
for each $w\in W'$.
\end{enumerate}

\begin{figure}[ht]
\centering
\includegraphics[width=9cm] {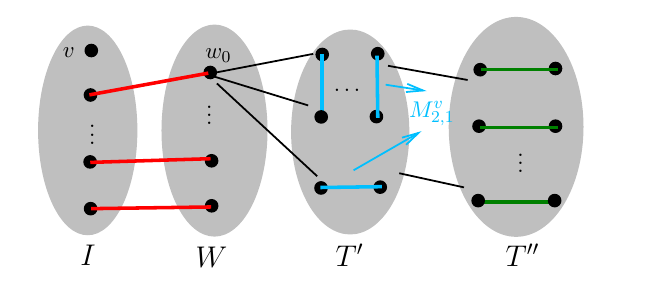}
\caption{$W'=\{w_0\}$ when $M_{2,2}^v=\emptyset$}
\label{figure4}
\end{figure}

In this subsection, we will show that
$W'=\{w_0\}$ and
$G[T''\cup T_{w_0}]$ is connected, as shown in Figure \ref{figure4}.
The notations $w_0$ and $W'$ will be applied in the remainder of this section.
Let $\setx=
\{T_{w_0}\cup T''\}\cup \{T_w: w\in W'\setminus \{w_0\}\}$.

\begin{lemma}\label{le3-2-01}
Let $G\in \setg$ and $w_0$ be the vertex in $W'$ introduced
in the beginning of Subsection 3.2.
Then, $T_w$ is a clique for each
$w\in W'\setminus \{w_0\}$.
\end{lemma}

\myproof
If $T_{w}$ is not a clique
for some $w\in W'\setminus \{w_0\}$,
by the assumption of $w_0$,
$T_{w_0}$ is not a clique,
implying that $\al(T_{w})\ge 2$
for at least two vertices $w$ in
$W$, and thus
$\al(T_{w})=2$
for at least two vertices $w$ in
$W$ by Lemma~\ref{basic-2} (iii),
 contradicting  Lemma~\ref{basic-4}.
Thus, the result holds.
\proofend

\begin{lemma}\label{le3-2-0}
	Let $G=(V,E)\in \setg$ with $M_{2,2}^v=\emptyset$.
	Assume that
	$w_0\in W'$ is the vertex
	in the beginning of Subsection 3.2.
	Then $G[T_{w_0}\cup T'']$ is connected.
\end{lemma}

\myproof Suppose that $G[T_{w_0}\cup T'']$ is not connected.
By definition, $G[T_{w_0}]$ is  connected, implying that
$T''\ne \emptyset$.
By Lemma~\ref{basic},
$T''$ is a clique,
implying that
$E(T_{w_0}, T'')=\emptyset$.
By Lemma~\ref{basic-3} (iii),
$T_{w_0}$ is a clique.
Thus,  Lemma~\ref{le3-2-01}
implies that
$T_w$ is a clique for each
$w\in W'$.

Since $E(T_{w_0}, T'')=\emptyset$,
we have $E(T_{w}, T'')=\emptyset$
for each $w\in W$ by the assumption
of $w_0$. By Lemma~\ref{basic-3} (ii),
$E(T_{w}, T_{w'})=\emptyset$
for each pair of $w,w'\in W'$.
It follows that
$|\partial(V\setminus I)|
=|\partial(T'')|+\sum\limits_{w\in W}|\partial(T_w)|$.

For each $w\in W'$,
since $T_w$ is a clique, we have
$3\le |T_w|\le  r-1$,
and thus
$|\partial(T_w)|=|T_w|(r-|T_w|+1)\ge 2r-2$.
By Lemma~\ref{basic} (ii),
$T''$ is a clique,
implying that $|T''|\le r$ by the fact
that $G$ is $r$-regular and $|T''|$ is even.
Thus, $|\partial(T'')|\ge
|T''|(r-|T''|+1)\ge r$.

Since $|\partial(I)|=|\partial(V\setminus I)|$,  we have
\equ{}
{
	r\al(G)=|\partial(I)|=
	|\partial(T'')|+\sum\limits_{w\in W}|\partial(T_w)|
	\ge r+|W'|(2r-2)+(|W|-|W'|)r
	=r\al(G)+(r-2)|W'|,
}
a contradiction. Thus, the result holds.
\proofend

\begin{corollary}\label{sec3-c2}
		Let $G\in \setg$ with $M_{2,2}^v=\emptyset$.
		Then, $G[X]$ is connected for each $X\in \setx$.
\end{corollary}

\myproof The result follows directly from
Lemmas~\ref{le3-2-01} and~\ref{le3-2-0}.
\proofend

\begin{lemma}\label{le3-2-1}
	Let $G\in \setg$ with $M_{2,2}^v=\emptyset$.
	Assume that
	$w_0$  is the vertex in $W'$
	in the beginning of Subsection 3.2
	 and
	$G[T_{w_0}\cup T'']$ is connected.
	For any $X\in \setx$,
	\vspace{-3 mm}
	\begin{enumerate}
		[itemsep=-4mm, parsep=0.5cm]
		\item 	if $|X|\ge r+1$,
		then $|\partial(X)|\ge r+1$;
		
		\item 	if $|X|\le r-1$,
		then either
		$|\partial(X)|\ge  2r+1$ or
		$G[X]$ is isomorphic to
		some  graph in $\{K_3, K_{r-1}, K_{r-1}-e\}$, where $r=6$
		when $G[X]\cong K_3$; and
	
		\item if $G[X]$ is
		not $2$-connected,
		then $|\partial(X)|\ge 2r-1$.
	\end{enumerate}
\end{lemma}

\myproof (i). Since $N_G(x)\cap I\ne \emptyset$ for any vertex $x$ of $X$,
$|X|\ge r+1$ implies that
$|\partial(X)|\ge |X|\ge r+1$.

(ii).
Let $X\in \setx$ with $|X|\le r-1$.
By the definition,  $|X|$ is odd.
Observe that
\equ{eq31}
{
	|\partial(X)|=r|X|-2|E(G[X])|\ge r|X|-|X|(|X|-1)=|X|(r-|X|+1).
}
Suppose that $|\partial(X)|\le 2r$.
Then, (\ref{eq31}) implies that
$2r\ge |X|(r-|X|+1)$.
It follows that $r\leq |X|+1+\frac{2}{|X|-2}$ where $3\le |X|\le r-1$ and $r\geq 6$.
Hence we have $r=6$ when $|X|=3$, and $|X|=r-1$ otherwise.

For the case $|X|=3$ and $r=6$, 
if $X$ is not a clique, then
the inequality of
(\ref{eq31}) is strict,
implying that 
$2\times 6>3\times (6-3+1)$,
a contradiction. Thus, $G[X]\cong K_3$ in this case.


For the case $|X|=r-1$,
if $|E(G[X])|\le {r-1\choose 2}-2$,
then
$$
2\times r\ge |\partial(X)|\ge r|X|-2\left ({|X|\choose 2}-2\right)
=r(r-1)-((r-1)(r-2)-4),
$$
a contradiction.
Thus, $|E(G[X])|\ge 
{r-1\choose 2}-1$, implying that
$G[X]\cong K_{r-1}$ or
$G[X]\cong K_{r-1}-e$.
The result holds.

(iii).
Let $X\in \setx$.
By Corollary~\ref{sec3-c2},
$G[X]$ is connected.
Assume that
$G[X]$ is not
$2$-connected.
By (ii), $|\partial(X)|\ge  2r+1>  2r-1$ when $|X|\le r-1$, then the result holds.

Now consider the case $|X|\ge r$.
Since $G[X]$ is not $2$-connected,
$X$ is not a clique, and thus
$X=T_{w_0}\cup T''$.
By Lemma~\ref{basic-3} (iv),
there exists a cut-vertex $u$ in $G[X]$ such that
$G[X]-u$ has exactly two components
whose vertex sets are
 $O_1$ and $O_2$
with the properties that  both $O_1\cup \{u\}$ and
$O_2$ are cliques of $G$.
If $|O_2|\ge r$, then for any vertex
$u'\in N_G(u)\cap O_2$,
$N_G(u')\cap I\ne \emptyset$,
implying that
$d_G(u')\ge |O_2|-1+2\ge r+1$, a contradiction.
Similarly, $|O_1|\le r-1$.
Observe that
\equ{le3-2-1-e2}
{
	|\partial(X)|\ge
	|O_1|(r-|O_1|)+1+|O_2|(r-|O_2|)
	\ge 2r-1.
}
The result holds.
\proofend

\begin{lemma}\label{le3-2-2}
	Let $G=(V,E)\in \setg$ with $M_{2,2}^v=\emptyset$.
	Assume that
	$E(T_{w_0}\cup T'', T_w)=\emptyset$ for each $w\in W'\setminus \{w_0\}$.
	If $G[T_{w_0}\cup T'']$ is $2$-connected,
	 then $|W'|\le 2$;
	 otherwise,
	$|W'|=1$.
\end{lemma}

\myproof
For any $w\in W'\setminus \{w_0\}$,
by Lemma~\ref{le3-2-01},
$T_w$ is a clique.
It follows that
$3\le |T_w|\le r-1$ and
$|\partial(T_w)|=|T_w|(r-|T_w|+1)\ge 2(r-1)$, where the equality holds
whenever $|T_w|=r-1$.
Note that $|\partial(I)|=|\partial(V\setminus I)|$.
For any $w_1,w_2\in W'\setminus \{w_0\}$,
$E(T_{w_1}, T_{w_2})=\emptyset$
by Lemma~\ref{basic-3} (ii).

Let $X_0=T_{w_0}\cup T''$.
Since 	$E(T_{w_0}\cup T'', T_w)=\emptyset$ for each $w\in W'\setminus \{w_0\}$, we have
\eqn{le3-2-2-e1}
{
r|I|&=&|\partial(I)|
=|E(I,X_0)|+
	\sum_{w\in W\setminus \{w_0\}}|E(I, T_w)|
\nonumber \\
&=&
|\partial(X_0)|+\sum_{w\in W'\setminus \{w_0\}}|\partial(T_w)|
+(|W|-|W'|)r
\nonumber \\
&\ge &
|\partial(X_0)|+(|W'|-1)(2r-2)+(|W|-|W'|)r
\nonumber \\
&= &
r(|W|+1)+(r-2)|W'|-(3r-2)+
|\partial(X_0)|.
}
Since $|I|=|W|+1$ by Corollary~\ref{nomatch-c1},
it follows from (\ref{le3-2-2-e1}) that
$|W'|\le \frac 1{r-2}
(3r-2-|\partial(X_0)|)$.

By Lemma~\ref{le3-2-1},
$|\partial(X_0)|\ge r+1$
if $G[X_0]$ is $2$-connected,
and  $|\partial(X_0)|\ge 2r-1$
otherwise.
Thus,
$|W'|\le \frac 1{r-2}
(3r-2-|\partial(X_0)|)$
implies that $|W'|\le 2$ if $G[X_0]$
is $2$-connected, and
$|W'|=1$ otherwise.
The result holds.
\proofend

\begin{lemma}\label{Nbasic-3}
Let $G=(V,E)\in \setg$ with $M_{2,2}^v=\emptyset$.
	If $G[T_{w_0}\cup T'']$ is not $2$-connected,
then $V=I\cup W\cup T_{w_0}\cup T''$,
i.e., $W'=\{w_0\}$.
\end{lemma}

\myproof
Let $X_0=T_{w_0}\cup T''$.
By Lemma~\ref{le3-2-0},
$G[X_0]$ is connected.
Now assume that $G[X_0]$
is not $2$-connected
and $W'\ne  \{w_0\}$.
By Lemma~\ref{le3-2-01},
$T_w$ is a clique for each $w\in W'\setminus\{w_0\}$.

Since $G[X_0]$ is not $2$-connected,
	by Lemma~\ref{basic-3}  (iv),
	$G[X_0]$ has a cut-vertex $u$  such that $G[X_0]-u$ has exactly two components
	whose vertex sets are
	$O_1$ and $O_2$ with the properties that both $O_1\cup\{u\}$ and
  	$O_2$ are cliques of $G$ and $G[X_0]-x$  has perfect matchings for each $x\in O_1\cup O_2$.
  We claim that
  	$E_{G}(O_1\cup O_2, T_w)=\emptyset$ for
  	each $w\in W\setminus \{w_0\}$.
  	Otherwise, there exists
  	$xy\in E_{G}(O_1\cup O_2, T_w)$,
  	where $x\in O_1\cup O_2$
  	and $y\in T_w$.
  	Then,  as $T_w$ is an odd clique
  	and $G[X_0]-x$  has perfect matchings,
  	$G[X_0\cup T_w]$ has perfect matchings,
  	 contradicting
  	Lemma~\ref{basic-3}  (i) .
	
  For each $w\in W\setminus\{w_0\}$,
 since $E_{G}(O_1\cup O_2, T_w)=\emptyset$, we have
$E_G(u,T_w)=E_G(X_0, T_w)$.
Since $G[X_0]$ is not $2$-connected,
by Lemma~\ref{le3-2-2},
there exists $w^*\in W'\setminus \{w_0\}$ such that
$E(X_0, T_{w^*})\ne \emptyset$.
It follows that both $|O_1|$ and $|O_2|$ are odd by the fact that  $E_G(u,T_{w^*})=E_G(X_0, T_{w^*})\neq\emptyset$. Otherwise, both $G[X_0]-u$ and $G[T_{w^*}]+u$ have perfect matchings,
implying that $G[X_0\cup T_{w^*}]$ has a perfect matching,  contradicting  Lemma~\ref{basic-3}  (i).

We are now going to apply
Proposition~\ref{move-T}.
Let $y_i\in O_i$ for $i\in\{1,2\}$.
Since $|O_i|$ is odd, $G[O_i]-y_i$ has perfect matchings for both $i\in\{1,2\}$.
Clearly,
$E_G(X_0, T_w)=E_G(u,T_w)$ implies that $N_{G}(y_i)\cap (W\setminus \{w_0\})=\emptyset$ for both $i\in\{1,2\}$.
It follows that
 $I':=\{y_1,y_2\}\cup (W\setminus \{w_0\})$ is an independent set of $G$
with $I'\cap I=\emptyset$.
By Corollary~\ref{nomatch-c1},
$|I'|=|W|+1=|I|=\al(G)$.

By the selection of $y_1$ and $y_2$,
$|N_G(z)\cap \{y_1,y_2\}|\le 1$
for each $z\in X_0\setminus
\{y_1, y_2, u\}$.
Hence, for each
$h\in V\setminus (I\cup I'\cup \{u\})$,
$|N_G(h)\cap I'|\le 1$ by Lemma \ref{basic-2}(ii) and the assumption of $T''$.
It is also clear that $G-(I\cup I'\cup \{u\})$ has a perfect matching
which is the union of
$\matt\setminus E(G[X_0])$
and two perfect matchings
in both $G[O_1]-y_1$ and
$G[O_2]-y_2$.

Let $X=T_w$ for some $w$ in
$W\setminus \{w_0\}$.
Obviously, $|X\cap I'|=1$.
For each
$w'\in W'\setminus \{w_0,w\}$,
since both $T_w$ and
$T_{w'}$ are cliques,
$E_G(X,T_{w'})=\emptyset$
by  Lemma~\ref{basic-3} (ii).
Since $X=T_w$ and
$E_G(T_w, X_0)
=E_G(T_w,u)$,
we have
$N_G(X)\subseteq \{u\}\cup X\cup I$.
However,
by Proposition~\ref{move-T},
such a set $X$ does not exist,
a contradiction.
Hence the result holds.
\proofend

\begin{lemma}\label{no2-conX}
	Let $G=(V,E)\in \setg$ with $M_{2,2}^v=\emptyset$.
	If $G[T_{w_0}\cup T'']$ is $2$-connected,
	then $V=I\cup W\cup T_{w_0}\cup T''$, i.e., $W'=\{w_0\}$,
	and $G[T_{w_0}\cup T'']$ is isomorphic to some graph in
	$\{K_3, K_{r-1}, K_{r-1}-e\}$,
	where $r=6$ if
	$G[T_{w_0}\cup T'']\cong K_3$.
\end{lemma}

\myproof Let $X_0=T_{w_0}\cup T''$.
Assume that $G[X_0]$ is $2$-connected.
By Lemma~\ref{le3-2-01},
$T_w$ is a clique for each $w\in W'\setminus\{w_0\}$.

\inclaim
$E_G(X_0, T_w)= \emptyset$ for each $w\in W\setminus\{w_0\}$.

Suppose that $E_G(X_0,T_w)\ne \emptyset$ for some $w\in W\setminus\{w_0\}$.
By Lemma~\ref{basic-3} (iii),
 $\al(G[X_0])\le 2$.
Since
$G[X_0]$ is $2$-connected,
by Corollary~\ref{equi-co},
$G[X_0]$ is factor-critical.
Then, for any edge $u_0u_1\in E_G(X_0,T_w)$,
where
$u_0\in X_0$ and $u_1\in T_w$,
both $G[X_0]-u_0$ and $G[T_w]-u_1$ have perfect matchings as $T_w$ is an odd clique,
implying that
$G[X_0\cup T_w]$ has a perfect matching,
 contradicting
Lemma~\ref{basic-3}  (i).
Thus, the claim holds.

By Claim 1 and Lemma~\ref{le3-2-2},
$|W'|\le 2$.
Thus, either $W'=\{w_0\}$ or
$W'=\{w_0,w_1\}$.
If $W'=\{w_0,w_1\}$,
let $X_1=T_{w_1}$.
Clearly,
$X_1$ is a clique and $|X_1|\leq r-1$. Then $\setx=\{X_0\}$ when $|W'|=1$,
and $\setx=\{X_0, X_1\}$ otherwise.

\inclaim $|\partial(X_0)|\le 2r$,
and if $\setx=\{X_0, X_1\}$, then
$|\partial(X_0)|+|\partial(X_1)|= 3r$.

If $X_1$ exists, then $X_1$ is a clique
and
$E(X_0, X_1)=\emptyset$ by Claim 1.
Thus,
\equ{eq32}
{
	r\al(G)
	=r\cdot |I|
	=|\partial(V\setminus I)|
	=|\partial(X_0)|+\sum_{X\in \setx\setminus \{X_0\}}|\partial(X)|
	+r\cdot (\al(G)-1-|W'|).
}
By (\ref{eq32}),
if $|\setx|=1$, then (\ref{eq32}) implies  that
$|\partial(X_0)|=2r$.
If $\setx=\{X_0, X_1\}$,
by (\ref{eq32}),
$|\partial(X_0)|+|\partial(X_1)|=3r$.
By  Lemma~\ref{le3-2-1} (ii),
it can be verified that
$|\partial(X_1)|\ge 2r-2$.
Thus, $|\partial(X_0)|\le r+2<2r$.
The claim holds.

\inclaim If $|X_0|\le r-1$, then $|W'|=1$
and the result holds.

Assume that $|X_0|\le r-1$.
By Lemma~\ref{le3-2-1} (ii),
either $|\partial(X_0)|\ge 2r+1$ or
$G[X_0]$ is isomorphic to
some  graph in $\{K_3, K_{r-1}, K_{r-1}-e\}$, where $r=6$ when
$G[X_0]\cong K_3$.
It follows that $|\partial(X_0)|\ge 2r-2$.
If $W'=\{w_0,w_1\}$ and
$\setx=\{X_0, X_1\}$, then
$X_1=T_{w_1}$ is a clique
by Lemma~\ref{le3-2-01},
implying that  $|X_1|\leq r-1$,
and thus $|\partial(X_1)|\ge 2r-2$ by Lemma~\ref{le3-2-1} (ii).
It follows that  $|\partial(X_0)|+|\partial(X_1)|>4r-4>3r$,
contradicting Claim 2.
 Hence $|W'|=1$ and $\setx=\{X_0\}$.
By Lemma~\ref{le3-2-1} (ii) and Claim 2 again, the claim holds.

Hence it remains to show that $|X_0|\ge r$ cannot happen.

\inclaim If $|X_0|\ge r+1$ and $\setx=\{X_0,X_1\}$, then $|X_1|=r-1$.

Assume that $|X_0|\ge r+1$.
By Lemma~\ref{le3-2-1} (i), $|\partial(X_0)|\ge r+1$.
Since $X_1$ is an odd clique and $G$ is $r$-regular,
 we have $|X_1|\leq r-1$.
Suppose that $|X_1| < r-1$.
Then Lemma~\ref{le3-2-1} (ii)
implies that
either
$|\partial(X_1)|\ge 2r+1$ or
$G[X_1]\cong K_3$ and $r=6$.
But, if
$G[X_1]\cong K_3$ and $r=6$,
$|\partial(X_1)|=3(r-2)=12=2r$.
Thus, we always have
$|\partial(X_1)|\ge 2r$.
Since $|\partial(X_0)|\ge r+1$
and Claim 1,
\equ{eq33}
{
r\al(G)=r\cdot |I|=|\partial(X_0)|+|\partial(X_1)|
+r(|W|-2)
\ge r+1+2r+r(\al(G)-3)>r\al(G),
}
a contradiction.  The claim  holds.

\inclaim If $|X_0|\ge r+1$, then
$N_G(z)\not\subseteq X_0$
for each $z\in I$.

Assume that $|X_0|\ge r+1$
and
$N_G(z)\subseteq X_0$ for some $z\in I$.  Since $|X_0|\ge r+1$ and $G$ is $r$-regular,
then
$N_G(z)\subseteq X_0\setminus \{u\}$ for some $u\in X_0$.
As $G[X_0]$ is 2-connected,
and $\al(G[X_0])\le 2$ by Lemma~\ref{basic-3} (iii),
$G[X_0]$ is factor-critical
by Corollary~\ref{equi-co}.
Thus,  $G[X_0]-u$
has  a minimal matching $M_0$
isolating $z$.
By Corollary \ref{Krr},
$G_0:=G-(V(M_0)\cup \{z\})$
is isomorphic to
either $K_{2t}$ or $K_{t,t}$
for some $t\ge 1$.
Note that $I_0:=I\setminus \{z\}$
and $W_0:=(W\setminus \{w_0\})\cup \{u\}$
are two independent sets in $G_0$
with  $|I_0|=|W_0|=\al(G)-1\ge 2$.
Thus, $G_0\cong K_{t,t}$ with $t\ge \al(G)-1$.
Since $u\notin N_{G}(z)$ and $N_{G}(u)\cap I\neq \emptyset$,
we have $I_0\subseteq N_G(u)$,
implying that
$|N_G(u)\cap I|\ge |I_0|=|I|-1$.
As $G[X_0]$ is 2-connected
and $u\in X_0$,
we have $|N_G(u)\cap X_0|\ge 2$.
Thus,
\equ{eq4}
{
	r=d_G(u)=|N_G(u)\cap X_0|+
	|N_G(u)\cap I|\ge 2+|I|-1=|I|+1,
}
implying that $|I|\le r-1$.

If $\setx=\{X_0,X_1\}$, then $X_1\subseteq V(G_0)$.
By Claim 4, $G[X_1]\cong K_{r-1}$,  contradicting the fact that $G_0$ is bipartite.
Hence $\setx=\{X_0\}$ and
$V=I\cup W\cup X_0$.
Note that
$\al\geq 3$ and $|W|=\al-1$ imply that $W\setminus \{w_0\}\neq \emptyset$.
For each $w\in W\setminus \{w_0\}$,
we have $E_G(T_w, X_0)=\emptyset$
by Claim 1, and
$N_G(w)\cap W=\emptyset$
by Corollary~\ref{nomatch-c1},
implying that
$N_G(w)\subseteq I$.
Hence $|I|\ge d_G(w)=r$,
a contradiction to (\ref{eq4}).
The claim  holds.

\inclaim $|X_0|\le r-1$.

Suppose that $|X_0|\ge r$. Since $|X_0|$ is odd, we have $|X_0|\ge r+1$.
By Lemma~\ref{basic} (i),
$|E_G(x,I)|\ge 1$ for each $x\in X_0$.
By Claim 2, $|\partial(X_0)|\le 2r$.
Thus,  $|X_0|\ge r+1$ implies that
$|E_G(x_0,I)|=1$ for some $x_0\in X_0$.
Let $I\cap N_G(x_0)=\{z_0\}$.
By Claim 5,
$N_G(z_0)\not\subseteq X_0$.
Choose a vertex $y_0$ in
$N_G(z_0)\setminus X_0$.

Since $G[X_0]$ is factor-critical by
Corollary~\ref{equi-co},
$G[X_0]-x_0$ has a perfect
matching $M'$. Let $M_0$ be
a minimum subset of  $M'$ such that
$N_{G}(x_0)\cap X_0\subseteq V(M_0)$.
Then $M'':=\{z_0y_0\}  \cup M_0$ is a minimal matching such that $x_0$ is isolated.
By Corollary \ref{Krr},
$G_1:=G-(V(M'')\cup \{x_0\})$
is isomorphic to either
$K_{2t}$ or $K_{t,t}$
for some $t\ge 1$.
Observe that both $I_1:=I\setminus \{z_0\}$ and $W_1:=W\setminus \{w_0,y_0\}$ are
independent sets of $G_1$.
As $|I_1|= \al(G)-1\ge 2$,
we have
$G_1\cong K_{t,t}$, where $t\ge \al(G)-1$.
Assume that $(A,B)$ is a bipartition of $G_1$ with $I_1\subseteq A$.
Clearly, $|A|=|B|$.

We claim that $\setx=\{X_0\}$.
If
$\setx=\{X_0, X_1\}$,  then $X_1=T_{w_1}$ is a clique
by Lemma~\ref{le3-2-01},
and $|X_1|=r-1$ by Claim 4.
As
$X_1\setminus \{y_0\}\subseteq V(G_1)$,
$G_1$ contains a clique of size  at least $r-2\ge 4$,
contradicting the fact that $G_1$ is bipartite.
Hence $\setx=\{X_0\}$.
Thus,  $N_{G}(z_0)\nsubseteq X_0$
imply that $N_{G}(z_0)\setminus X_0 \subseteq W\setminus\{w_0\}$ and $y_0\in W\setminus\{w_0\}$.

Note that $V(G_1)=I_1\cup W_1\cup
(X_0\setminus (V(M_0)\cup \{x_0\}))
$.
By the assumption of $M_0$, $M'\setminus M_0$ is a perfect matching of $G[X_0]-(V(M_0)\cup \{x_0\})$.
Since $G_1$ is bipartite,
for each edge $u_1u_2$
in $M'\setminus M_0$,
we have $u_i\in A$ and $u_{3-i}\in B$
for some $i\in \{1,2\}$.
As $I_1\subseteq A$, and $|I_1|=|W_1|+2$,
we have $|A|\ne |B|$,
contradicting the fact that $|A|=|B|$,
no matter whether
$W_1\subseteq A$ or
$W_1\subseteq B$.
Thus Claim 6 holds.

By Claims 6 and 3, the result holds.
\proofend

\section{Proof of  $M_{2,2}^v\ne \emptyset$
	\label{sec4}
}

\begin{figure}[ht]
\centering
\includegraphics[width=7cm] {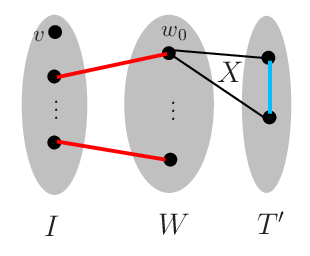}
\caption{$T_{w_0}$ is a clique
	of size $3$ and $N_{G}(T_{w_0})\subseteq T_{w_0}\cup I$}
\label{figure5}
\end{figure}

In this section, we will show that  $M_{2,2}^v\neq \emptyset$.
We first prove that if $M_{2,2}^v= \emptyset$, then $r=6$
and $V=I\cup W\cup T_{w_0}$, where
$w_0$ is the vertex introduced
in the beginning of Subsection 3.2, and
$T_{w_0}$ is a clique with
$|T_{w_0}|=3$ and $N_{G}(T_{w_0})\subseteq T_{w_0}\cup I$, as shown in Figure \ref{figure5}. Next we show that such a structure
cannot exist.

\begin{lemma}\label{Nbasic-4}
Let $G=(V,E)\in \setg$ with $M_{2,2}^v=\emptyset$.
Then,  $r=6$
and $V=I\cup W\cup X$,
where $X\subseteq V\setminus I$
is a clique with $|X|=3$,
$|X\cap W|=1$
and $N_G(X)\subseteq X\cup I$,
as shown in Figure \ref{figure5}.
\end{lemma}

\myproof
Let $w_0$ be the vertex introduced
in the beginning of Subsection 3.2. 
By Lemmas~\ref{Nbasic-3} and~\ref{no2-conX},
$W'=\{w_0\}$ and
$V=I\cup W\cup T_{w_0}\cup T''$.
Let $X=T_{w_0}\cup T''$.
By Lemma~\ref{le3-2-0},
$G[X]$ is connected.
By the definition of $T''$ and Lemma~\ref{basic-2} (ii), $E(X, W\setminus \{{w_0}\})=\emptyset$.

It is known that $G[X]$ is connected.
Since $W'=\{w_0\}$ and $E(X, W\setminus \{{w_0}\})=\emptyset$, we have
 \equ{eq11}
{
	r\al(G)=|\partial(I)|=
	|\partial(X\cup (W\setminus \{w_0\}))|=
	|E_G(X,I)|+r|W\setminus \{{w_0}\}|
	= |\partial(X)|
	+r(\al(G)-2),
}
implying that $|\partial(X)|=2r$.

\inclaim
$G[X]$ is 2-connected.

Suppose $G[X]$ is not $2$-connected.
By Lemma~\ref{basic-3} (iv),
$G[X]$ has a
cut-vertex $u$ such that
$G[X]-u$ has exactly two components
whose vertex sets are
 $O_1$ and $O_2$
with the property that
both $O_1\cup \{u\}$ and
$O_2$ are cliques of $G$.
Since $G$ is connected and $r$-regular,
$|O_i|\le r$ for both $i\in\{1,2\}$.
If $|O_1|=r$,  then
the facts that $O_1$ is a clique,
and $G$ is $r$-regular
imply that $N_G(z)\cap I=\emptyset$ for some $z\in N_G(u)\cap O_1$,
contradicting Lemma~\ref{basic} (i). Thus, $|O_1|\le r-1$. Similarly,
$|O_2|\le r-1$.

\begin{figure}[ht]
	\centering
	\includegraphics[width=8 cm]
	{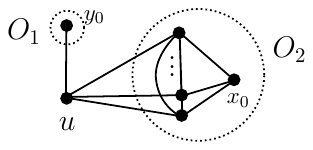}
	
	\caption{$O_1=\{y_0\}$, $O_2$ is a clique of size $r-1$ and $N_G(u)\cap O_2= O_2\setminus \{x_0\}$}
	
	\label{Nbasic-4-f1}
\end{figure}

Since $O_1\subseteq N_G(u)$,
$N_G(u)\cap I\ne \emptyset$
and $N_G(u)\cap O_2\ne \emptyset$,
we have $|O_1|\le r-2$.
Since $O_i$ is a clique
and $N_G(O_i)\subseteq O_i\cup \{u\}\cup I$ for 
both $i\in\{1,2\}$,
 \equ{eq12}
{
|E_G(O_1, I)|= |O_1|(r-|O_1|),
\quad
\mbox{and}\quad
|E_G(O_2, I)|\ge |O_2|(r-|O_2|).
}
Clearly, $|O_i|(r-|O_i|)=r-1$ when
$|O_i|\in \{1,r-1\}$,
and $|O_i|(r-|O_i|)\ge 2r-4$
when $2\le |O_i|\le r-2$.
Thus, if $2\leq|O_1|\leq r-2$, then
$|O_1|(r-|O_1|)\ge 2r-4$,
implying that
 \equ{eq13}
{
|\partial(X)|=|E_G(u,I)|+
|E_G(O_1, I)|+|E_G(O_2, I)|
\ge 1+(2r-4)+(r-1)=3r-4>2r,
}
contradicting (\ref{eq11}).
If $|O_1|=1$ and
$2\le |O_2|\leq r-2$, then
 \equ{eq14}
{
	|\partial(X)|=|E_G(u,I)|+
	|E_G(O_1, I)|+|E_G(O_2, I)|
	\ge 1+r-1+2r-4=3r-4>2r,
}
contradicting (\ref{eq11}).
If $|O_1|=|O_2|=1$, then
 \equ{eq15}
{
|\partial(X)|=|E_G(u,I)|+
|E_G(O_1, I)|+|E_G(O_2, I)|
=(r-2)+(r-1)+(r-1)=3r-4>2r,
}
contradicting (\ref{eq11}).
Finally, for the case that $|O_1|=1$ and $|O_2|=r-1$,
$N_{G}(u)\cap I\neq \emptyset$ and $O_1\subseteq N_{G}(u)$ imply that $|N_{G}(u)\cap O_2|\leq r-2$.
As $O_2$ is a clique of size $r-1$,
we have
$|E_G(O_2, I)|=
(2(r-1)-|N_G(u)\cap O_2|)$,
implying that
\eqn{eq2}
{
	|\partial(X)|
	&=&|E_G(u,I)|+
	|E_G(O_1, I)|+|E_G(O_2, I)|
	\nonumber \\
	&=&(r-1-|N_G(u)\cap O_2|)+(r-1)
	+(2(r-1)-|N_G(u)\cap O_2|)
	\nonumber \\
	&= &4(r-1)-2|N_G(u)\cap O_2|\geq 2r,
}
where the equality holds if and only
if $|N_G(u)\cap O_2|=r-2$.
If $|O_1|=1$, $|O_2|=r-1$
and
$|N_G(u)\cap O_2|=r-2$,
$G[X]$ is
the graph shown in Figure~\ref{Nbasic-4-f1}.

Now we are going to apply Proposition~\ref{move-T} to show that
the above conclusion leads to a contradiction.
Let $I'=\{x_0,y_0\}\cup (W\setminus \{w_0\})$,
where $y_0$ is the only vertex in $O_1$ and $x_0$ is the only vertex in
$O_2\setminus N_G(u)$, see Figure~\ref{Nbasic-4-f1}.
Clearly, $I'\cap I=\emptyset$,
and $I'$ is an independent set of $G$.
Observe that
$G-(I\cup I'\cup \{u\})=G[O_2\setminus \{x_0\}]$
has a perfect matching, since it is
actually the graph
$G[O_2\setminus \{x_0\}]\cong K_{r-2}$.
As
$V_0:=V\setminus
(I\cup I'\cup \{u\})=O_2\setminus \{x_0\}$,
for any vertex $u'\in V_0$,
we have $N_G(u')\cap I'=\{x_0\}$.
Now let $X':=O_2$.
Observe that  $X'$ is an odd clique
of $G-(I\cup \{u\})$
with the property that
$|X'|=r-1>3$, $X'\cap I'=\{x_0\}$,
and $N_G(X')\subseteq
\{u\}\cup X'\cup I$.
By Proposition~\ref{move-T},
such a clique $X'$ of $G-(I\cup\{u\})$ should not exist, a contradiction.
Thus, Claim 1 holds.

Now we are going to complete the proof.
By Claim 1, $G[X]$ is $2$-connected.
Then, by Lemma~\ref{no2-conX},
$V=I\cup W\cup X$,
$G[X]\cong K_3$ with $r=6$,
or $G[X]\cong K_{r-1}$ or
$G[X]\cong K_{r-1}-e$.
If $G[X]\cong K_{r-1}$, we have $|\partial(X)|=2r-2$, contradicting
(\ref{eq11}).
In the following, we assume that
$G[X]\cong K_{r-1}-e$,
and we will show that
this conclusion
leads to a contradiction.

\inclaim 
There is no matching
$M_0$ of $G$
with $M_0\subseteq
E_G(I,X)$ and $|M_0|=4$.

Suppose that
$G$ has a matching $M_0\subseteq E_G(I, X)$ with $|M_0|=4$.
Then, $G[X\setminus V(M_0)]-x'$ has a perfect matching $M_1$ for
some $x'\in X\setminus V(M_0)$ due to the fact that $G[X]\cong K_{r-1}-e$.
Note that the vertex set of $G-V(M_0\cup M_1)$ can be partitioned into independent sets
$I_1:=I\setminus V(M_0)$
and $W_1:=(W\setminus X)\cup \{x'\}$.
Observe that $|I_1|=\al(G)-4$
while $|W_1|=\al(G)-1$,
contradicting Proposition~\ref{bipartite}.
Thus, Claim 2 holds.

Let $x_1$ and $x_2$ be
the two non-adjacent vertices in $G[X]$.
Then, for any $x\in X$,
$|N_G(x)\cap I|=3$ if $x\in \{x_1,x_2\}$,
and $|N_G(x)\cap I|=2$  otherwise.
By Claim 2,
Hall's Theorem \cite{hall1987}
implies the existence of 
	a subset $I_0$ of 
$I$ with $|I_0|=3$ such that
$N_G(x)\cap I\subseteq I_0$
for each $x\in X$.
Clearly,
$N_G(x_i)\cap I=I_0$ for both $i\in\{1,2\}$.

\inclaim
There is no matching
$M_0$ with $M_0\subseteq E_G(I_0, X\setminus \{x_1,x_2\})$ and
$I_0\subseteq V(M_0)$.

Suppose that $M_0\subseteq E_G(I_0, X\setminus \{x_1,x_2\})$
is a matching such that
$I_0\subseteq V(M_0)$.
Then,  $|M_0|=3$ and
the graph
$H:=G[X]-(V(M_0)\cap X)$
is isomorphic to $K_{r-4}-e$
and the two non-adjacent vertices in $H$ are $x_1$ and $x_2$.
Thus, $H-\{x_1,x_2\}$
has a perfect matching $M'$.
Note that $M'=\emptyset$ when $r=6$.
Then, $x_1$ is isolated by
the minimal matching $M'':=M_0\cup M'$.
By Corollary \ref{Krr},
$G-(V(M'')\cup \{x_1\})$
is isomorphic to $K_{2t}$ or $K_{t,t}$
for some $t\ge 1$.
However, the vertex set of
$G-(V(M'')\cup \{x_1\})$  can be partitioned into independent sets
$I\setminus V(M_0)$
and $(W\setminus X)\cup \{x_2\}$ whose sizes differ by 2, a contradiction.
Thus, Claim 3 holds.

\begin{figure}[ht]
	\centering
	\includegraphics[width=7 cm]{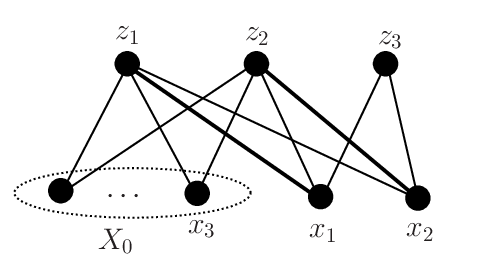}
	
	\caption{$N_G(x)\cap I=\{z_1,z_2\}$
		for each $x\in X_0$, and
		$N_G(x_i)\cap I=\{z_1, z_2, z_3\}$ for both $i\in \{1,2\}$
	}
	
	\label{d-Wge1-f1}
\end{figure}

Claim 3 implies that
there exists $I':=\{z_1,z_2\}
\subseteq I_0$ such that
$N_G(x)\cap I=I'$
for each $x\in X_0:=X\setminus \{x_1,x_2\}$, as shown in Figure~\ref{d-Wge1-f1},
where $\{z_3\}=I_0\setminus I'$.
Let $x_3$ be a vertex in $X_0$.
Clearly, $G[X]-\{x_i:i\in [3]\}$ is isomorphic to $K_{r-4}$ and thus
has a perfect matching $N_1$.
Then, $N_2:=\{x_1z_1,x_2z_2\}\cup N_1$
is a minimal matching of $G$
isolating $x_3$.
By Corollary \ref{Krr},
$H_0:=G-(V(N_2)\cup \{x_3\})$
is isomorphic to $K_{2t}$ or $K_{t,t}$.
Observe that
$V(H_0)$ is partitioned into independent sets
$I_2:=I\setminus \{z_1,z_2\}$
and $W_2:=W\setminus X$. Thus,
$H_0\cong K_{t,t}$, where
$t=\al(G)-2$.
As $z_3\in I_2$,
we have $W_2\subseteq N_G(z_3)$.
Since $N_G(z_3)\cap X=\{x_1,x_2\}$,
 \equ{eq16}
{
r=d_G(z_3)=2+|W_2|=2+\al(G)-2
=\al(G),
}
implying that $\al(G)=r$.
Let $z\in I\setminus \{z_1,z_2,z_3\}=I\setminus I_0$.
Note that $N_G(z)\subseteq W\cup X$.
Since $N_G(x)\subseteq X\cup I_0$
for each $x\in X$,
we have
$N_G(z)\cap X=\emptyset$,
implying that
$N_G(z)\subseteq W\setminus X=W_2$.
It follows that
$\al-2=|W_2|\ge d_G(z)=r$,
implying that $\al\ge r+2$,
a contradiction to (\ref{eq16}).

Hence the result holds.
\proofend

Now we are going to establish the main
result in this section.

\begin{proposition}\label{M22}
For any $G=(V,E)\in \setg$,  we have $M_{2,2}^v\neq \emptyset$.
\end{proposition}

\myproof Suppose $M_{2,2}^v=\emptyset$.
By Lemma \ref{Nbasic-4},
we have $r=6$
and $V=I\cup W\cup X$,
where $X\subseteq V\setminus I$
is a clique with $|X|=3$,
$|X\cap W|=1$
and $N_G(X)\subseteq X\cup I$, as shown in Figure \ref{figure5}.

Let $X=\{w_0, x_1,y_1\}$, where $w_0\in W$.
Note that $N_G(X)\subseteq X\cup I$
by Lemma~\ref{basic-2} (ii).
Let $I_0=N_G(w_0)\cap I$.
Clearly, $|I_0|=r-2=4$
and $G_0:=G-\{w_0,x_1,y_1\}$
is bipartite.
If there exists a matching $M'$
in $G_0$
with $I_0\subseteq V(M')$ and
$|I_0|=|V(M')|/2$,
then $M'':=M'\cup \{x_1y_1\}$
is a minimal matching in $G$
isolating $w_0$.
By Corollary \ref{Krr},
$G-(V(M'')\cup \{w_0\})$
is isomorphic to $K_{2t}$ or $K_{t,t}$
for some $t\ge 1$.
However,  the vertex set of
$G-(V(M'')\cup \{w_0\})$
can be partitioned into two independent sets
$I\setminus V(M')$ and
$W\setminus (\{w_0\}\cup V(M'))$
whose  sizes differ by $2$,
a contradiction.

Thus, $G_0$ has no  matching $M'$
with $I_0\subseteq V(M')$ and
$|I_0|=|V(M')|/2$.
For each $z\in I_0$,
since
$N_G(z)\subseteq W\cup \{x_1,y_1\}$,
we have
$|N_G(z)\cap (W\setminus \{w_0\})|\ge r-3=3$.
Thus, $|W_0|\ge 3$, where
$W_0:=
N_{G_0}(I_0)\cap
(W\setminus \{w_0\})$.
Since $G_0$ has no  matching $M'$
with $I_0\subseteq V(M')$ and $|I_0|=4$, we have $|W_0|<|I_0|=4$,
i.e. $|W_0|=3$,
and
$N_G(z)=W_0\cup \{w_0,x_1,y_1\}$
for each $z\in I_0$.

Write $I_0=\{z_i:i\in [4]\}$ and
$W_0=\{w_1,w_2,w_3\}$.
Clearly, $w_0$ is isolated by the
matching
$M_0:=\{x_1z_1,y_1z_2, w_1z_3,w_2z_4\}$.
By Corollary \ref{Krr},
$G-(V(M_0)\cup \{w_0\})$
is isomorphic to $K_{2t}$ or $K_{t,t}$
for some $t\ge 1$.
Observe that
the vertex set of $G-(V(M_0)\cup \{w_0\})$
can be partitioned into two independent
sets $I_1:=I\setminus I_0$ and
$W_1:=W\setminus \{w_0,w_1,w_2\}$.
Thus, $t=|I_1|=\al(G)-4$, and so $\al(G)=4+t$.

Note that $w_3\in W_1$ and
$N_G(w_3)\subseteq I$.
Since $I_0\subseteq N_G(w_3)$
and $|I_0|=4$,
we have $|N_G(w_3)\cap I_1|=|I_1|=2$,
implying that $t=|I_1|=2$
and $\al(G)=6$.
As $|W|=\al(G)-1=5$,
there is a vertex $w'\in W\setminus
\{w_0,w_1,w_2,w_3\}$.
Since $N_G(I_0)\cap W=\{w_i: i\in \{0,1,2,3\}\}$,
we have $N_G(w')\subseteq I\setminus I_0=I_1$,
implying that $|I_1|\ge d_G(w')=6$,
contradicting $|I_1|=2$.
Hence the result holds.
\proofend

\section{Proof of Theorem~\ref{main1} \label{sec5}
}

By Proposition \ref{M22} and Lemma~\ref{basic-2} (iv), $M^v_{2,2}$ contains
exactly one edge, say $uu'$.
We first deduce some
conclusions
on $u,u'$ and $T_w$'s for $w\in W'$.

\begin{lemma}\label{yT''}
Let $G\in \G$ and
$M^v_{2,2}=\{uu'\}$ with
$|N_G(u)\cap W|\ge 2$.
Then,
\vspace{-3 mm}
\begin{enumerate}
	[itemsep=-4mm, parsep=0.5cm]
\item $N_G(u')\cap W=\emptyset$;

\item $T''\cup \{u'\}$ is a clique;

\item for any $w\in W'$,
$T_{w}$ is a clique and
$N_{G}(u')\cap T_{w} =\emptyset$;
 and

\item $N_G(T'')\subseteq
T''\cup I\cup
\{u,u'\}$.

\end{enumerate}
\end{lemma}

\myproof
Since $|N_G(u)\cap W|\ge 2$,
there are two vertices
$w_1$ and $w_2$ in
$N_G(u)\cap W$.

(i). If $w'\in N_G(u')\cap W$, then there is a path $w_iuu'w'$
for some $i\in \{1,2\}$ with
$w_i\ne w'$, a contradiction to
Lemma~\ref{basic-2} (i).
Thus, (i) holds.

(ii). By (i), $N_G(u')\cap W=\emptyset$.
Thus, $u'\in T'\setminus N_{G}(W)$,
implying that  $T''\cup \{u'\}$
is a clique by Lemma~\ref{basic} (ii).

(iii).  Recall that $W'$ denotes the set of vertices $w\in W$ for which $|T_w|\ge 3$.
	 We first prove the following claim.

\inclaim  For any $w\in W$,
$\al(G[T_{w}\cup \{u'\}])\le 2$.

Suppose that Claim 1 fails.
Then, $G[T_{w}\cup \{u'\}]$ has
an independent set $I_0$
with $|I_0|\ge 3$.
Since $N_G(u')\cap W=\emptyset$ by (i),
Lemma~\ref{basic-2} (ii)
implies that
$I_0\cup (W\setminus \{w\})$
is an independent set
of $G$ with $|I_0\cup (W\setminus \{w\})|=\al(G)+1$, a contradiction.
The claim holds.

\begin{figure}[ht]
	\centering
	\includegraphics[width=12 cm]
	{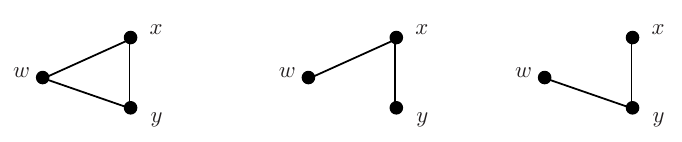}
	
(a) 	\hspace{3.6 cm}
(b) \hspace{3.7 cm}
(c)

	\caption{Three possible types of  edges $xy\in M^v_2\setminus \{uu'\}$ with $\{x,y\}\cap N_G(w)\ne \emptyset$}
	
	\label{figure8}
\end{figure}

Let $w\in W'$.
The set $T_w\setminus \{w\}$ can be partitioned into three subsets:

(i).   $S_1(w):=\{x,y:
xy\in M^v_2\setminus \{uu'\}, \{x,y\}\subseteq N_G(w)\}$,
as shown in Figure~\ref{figure8}
(a);

(ii).  $S_2(w):=\{x: xy\in
M^v_2\setminus \{uu'\}, \{x,y\}\cap N_G(w)=\{x\}\}$,
as shown in Figure~\ref{figure8}
(b);
and

(iii).  $S_3(w):=\{x: xy\in
M^v_2\setminus \{uu'\}, \{x,y\}\cap N_G(w)=\{y\}\}$,
as shown in Figure~\ref{figure8}
(c).

Since $|W\cap N_G(u)|\ge 2$,
there exists
$w_0\in (W\cap N_G(u))\setminus \{w\}$.
If there exists $x$ in $N_G(u')\cap (S_1(w)\cup S_3(w))$,
then $w_0uu'xyw$ is an $\matt$-augmenting path connecting $w_0$ and $w$,
where $y$ is the vertex with $xy\in \matt$,
contradicting
Lemma~\ref{basic-2} (i).
Hence $N_G(u')\cap (S_1(w)\cup S_3(w))=\emptyset$.

We claim that $S_2(w)=\emptyset$.
Otherwise, there exists $xy\in \matt$
with $xy\ne uu'$ and $\{x,y\}\cap N_G(w)=\{x\}$.
Then $y\in S_3(w)$.
As $N_G(u')\cap S_3(w)=\emptyset$,
$u'y\notin E(G)$,
implying that $\{w,u',y\}$ is an independent set,
contradicting Claim 1.

 By (i),  $w\notin N_{G}(u')$.
 Since $N_G(u')\cap (S_1(w)\cup S_3(w))=\emptyset$
 and $S_2(w)=\emptyset$,
 we have
 $N_G(u')\cap T_w=\emptyset$.
If $G[T_w]$ has
an independent set $I_0$ with $|I_0|=2$,
then $N_G(u')\cap T_w=\emptyset$
implies that
$I_0\cup \{u'\}$ is an independent set of size $3$, also contradicting Claim 1.
Hence $G[T_w]$  is a clique and
(iii) holds.

(iv). It suffices to show that
$E_G(T'',T_w)=\emptyset$ for each $w\in W'$.

Suppose that $yz\in E_G(T'',T_w)$,
where $w\in W'$, $y\in T_w$
and $z\in T''$.
Let $w_0\in N_G(u)\cap W\setminus \{w\}$.
By (ii) and (iii),  both
$T_w$ and $T''\cup \{u'\}$ are odd cliques. It follows that
$G[T''\cup T_w\cup \{w_0,u,u'\}]$
has a perfect matching $M_0$
with $w_0u, yz\in M_0$.
Let $M'=M_0\cup (M^{v}_2\setminus E(G[T''\cup T_w\cup \{w_0,u,u'\}]))$.
Then $V(G-V(M'))$
can be partitioned into two
independent sets $I$
and $W\setminus\{w_0, w\}$.
As $|W\setminus\{w_0, w\}|=|I|-3$,
we have $\al(G-V(M'))=|I|$
and $|V(G-V(M'))|=2|I|-3$.
However,  by
Proposition \ref{bipartite},
we have
$\al(G-V(M'))\le \frac{|V(G-(M'))|+1}{2}$,
i.e.,
$|I|\le \frac{2|I|-3+1}{2}=|I|-1$,
a contradiction.

Hence $E_G(T'',T_w)=\emptyset$ for each $w\in W'$
and (iv) holds.
\endproof

In order to prove Theorem~\ref{main1},
we are now going to give
a necessary and
sufficient condition
for a connected $r$-regular and equimatchable graph
to be  $F_r$.

\begin{lemma}\label{char-Fr}
	Let $G=(V,E)$ be a connected $r$-regular and equimatchable graph,
	where $r\ge 6$ is even.
	Then $G\cong F_r$ if and only if
	for some  vertex $u$,
	$V\setminus \{u\}$ can be partitioned into  independent sets $X$ and $Y$
	with $|X|=|Y|$.
\end{lemma}

\myproof The necessity is obvious.
We are now going to prove the sufficiency.

Assume that $u$ is a vertex in $V$
such that
$V\setminus \{u\}$ is partitioned
into two  independent sets $X$ and $Y$ with $|X|=|Y|$, and thus
$G-u$ is a bipartite graph
with a bipartition $(X,Y)$.
Let $X_0=N_G(u)\cap X$ and $Y_0=N_G(u)\cap Y$.
Then $|X_0|+|Y_0|=d_G(u)=r$.
We first prove the following claims.

\inclaim $|X_0|=|Y_0|=\frac r2$ and $|X|=|Y|\ge r$.

Since $G$ is $r$-regular, and both
$X$ and $Y$ are independent sets,
\equ{eq20}
{
	|E_G(u,X)|=
	r|X|-|E_G(X, Y)|
	=r|Y|-|E_G(X, Y)|=|E_G(u,Y)|,
}
implying that
$|X_0|=|E_G(u,X)|=|E_G(u,Y)|=
|Y_0|$.
Since $|X_0|+|Y_0|=r$, we have
$|X_0|=|Y_0|=\frac r2$.
For any $x\in X$, we have $N_G(x)\subseteq Y\cup \{u\}$,
implying that $r\le |Y|+1$, and so $|Y|\ge r-1$.
It follows that $|X|+|Y|\ge 2r-2>r=d_G(u)$,
and thus
$X\cup Y\not \subseteq N_G(u)$.
Assume that $y\in Y\setminus N_G(u)$.
As $Y$ is an independent set of $G$,
we have $N_G(y)\subseteq X$, implying that
$r=d_G(y)\le |X|=|Y|$.
Thus, Claim 1 holds.

\inclaim For any $x_1,x_2\in X_0$,
$G-\{u,x_1,x_2\}$ has no matching
$M_0$ with
$V(M_0)\cap Y=N_G(x_1)\setminus \{u\}$.
Similarly,
for any $y_1,y_2\in Y_0$,
$G-\{u,y_1,y_2\}$ has no matching
$M_0$ with
$V(M_0)\cap X=N_G(y_1)\setminus \{u\}$.

Suppose the claim fails, and that
$G-\{u,x_1,x_2\}$ has
a matching $M_0$ such that
$V(M_0)\cap Y=N_G(x_1)\setminus \{u\}$.
It follows that $|M_0|=d_G(x_1)-1=r-1$.
Then $M_1:=\{ux_2\}\cup M_0$
is a minimal matching
of $G$ isolating $x_1$.
By Corollary \ref{Krr},
$G-(V(M_1)\cup \{x_1\})$
is isomorphic to
$K_{2t}$ or $K_{t,t}$
for some positive integer $t$.
However,
the vertex set of
$G-(V(M_1)\cup \{x_1\})$
can be partitioned into
two independent sets
$X\setminus
(V(M_0)\cup \{x_1,x_2\})$
and $Y\setminus V(M_0)$
with the following property:
$$
|Y\setminus V(M_0)|-
|X\setminus
(V(M_0)\cup \{x_1,x_2\})|
=(|Y|-(r-1))-(|X|-(r+1))
=2,
$$
a contradiction
to the fact that
$G-(V(M_1)\cup \{x_1\})$
is isomorphic to
$K_{2t}$ or $K_{t,t}$ for some positive integer $t$.
Hence Claim 2 holds.

Let $X_1=(N_G(Y_0)\cap X)\setminus X_0$,
$Y_1=(N_G(X_0)\cap Y)\setminus Y_0$, $X_2=X\setminus (X_0\cup X_1)$
and $Y_2=Y\setminus (Y_0\cup Y_1)$.
By Claim 1,  $|X_0|=|Y_0|=\frac{r}{2}<r-1\le |X|=|Y|$. Thus $X_1\neq\emptyset$ and $Y_1\neq\emptyset$.

\begin{figure}[ht]
	\centering
\includegraphics[width=9 cm] {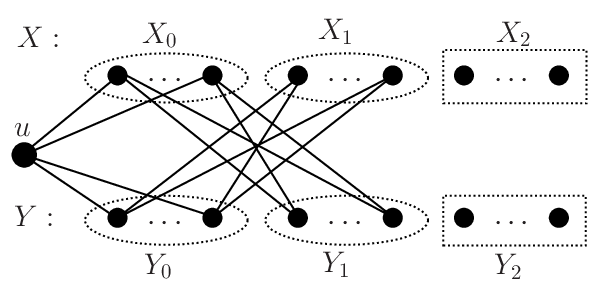}
	
	\caption{$X_1=(N_G(Y_0)\cap X)\setminus X_0$,
	$Y_1=(N_G(X_0)\cap Y)\setminus Y_0$, $X_2=X\setminus (X_0\cup X_1)$
		and $Y_2=Y\setminus (Y_0\cup Y_1)$}
	
	\label{char-Fr-f2}
\end{figure}

\inclaim
 Both $G[X_1\cup Y_0]$
and $G[Y_1\cup X_0]$
are complete bipartite graphs,
as shown in Figure~\ref{char-Fr-f2}.

Suppose that $G[Y_1\cup X_0]$
is not a complete bipartite graph.
Then,  there exists $y_0\in Y_1$
such that $X_0\not\subseteq N_G(y_0)$.
By the assumption of $Y_1$,
there exists $x_1\in X_0\cap
 N_G(y_0)$. Assume that
$x_2\in X_0\setminus N_G(y_0)$.

Let $Y':=N_G(x_1)\cap Y$
and $G_0:=G-\{u,x_1,x_2\}$.
Obviously, $y_0\in Y'$ and $N_{G_0}(y_0)=N_G(y_0)\setminus \{x_1\}$.
Since $\{u,x_2\}\cap N_G(y_0)=\emptyset$,
we have $d_{G_0}(y_0)=d_G(y_0)-1=r-1$.
Since $x_1\in X_0\subseteq N_G(u)$,
we have $|Y'|=d_G(x_1)-1=r-1$.
Thus, $|Y'\setminus Y_0|=(r-1)-\frac{r}{2}\ge 2$.
Let $y_1\in Y'\setminus (Y_0\cup \{y_0\})$.
Since $y_1\in Y\setminus Y_0=Y\setminus N_G(u)$,
we have $N_G(y_1)\subseteq X$,
implying that
$d_{G_0}(y_1)
=d_G(y_1)-|N_G(y_1)\cap \{x_1,x_2\}|\ge r-2$.

For any $y\in Y'$, we also have
$d_{G_0}(y)\ge d_G(y)-|N_G(y)\cap \{u,x_1,x_2\}|\ge r-3$.
We claim that for every subset $S\subseteq Y'$, $|N_{G_0}(S)|\ge |S|$.
It holds obviously
	for any $S\subseteq Y'$ with
$|S|\le r-3$.
Now assume that $r-2\le |S|\le r-1$.
Then $\{y_0, y_1\}\cap S\neq\emptyset$,
and $\{y_0, y_1\}\subseteq S$
if $|S|=r-1$.
Since $d_{G_0}(y_0)=r-1$ and $d_{G_0}(y_1)\ge r-2$, we can conclude that
$|N_{G_0}(S)|\ge |S|$.
By  Hall's Theorem \cite{hall1987},
$G_0$ has a matching $M'$ with
$V(M')\cap Y=Y'=N_G(x_1)\setminus \{u\}$,
contradicting
Claim 2.

Hence
$G[Y_1\cup X_0]$ is a complete bipartite graph.
Similarly,
$G[X_1\cup Y_0]$ is a complete bipartite graph.
Claim 3 holds.

\inclaim For any $x'\in X_0$,
$|N_G(Y')\cap X|\le r$,
	where $Y'=N_G(x')\cap Y$.
Similarly, for any $y'\in Y_0$,
$|N_G(X')\cap Y|\le r$,
where $X'=N_G(y')\cap X$.

Let $x',x\in X_0$ and $Y'=N_G(x')\cap Y$.
By Claim 2,
$G-\{u,x,x'\}$ does not have
a matching $M'$ such that
$V(M')\cap Y=Y'$.
Thus, according to Hall's Theorem \cite{hall1987},
there exists $S\subseteq Y'$ such that $|N_G(S)\cap (X\setminus \{x,x'\})|<|S|$.
Note that
$G-\{u,x,x'\}$  is a bipartite graph
with a bipartition $(Z,Y)$,
where $Z:=X\setminus \{x,x'\}$,

Suppose that
$|N_G(Y')\cap X|\ge r+1$.
Then $|N_G(Y')\cap Z|\ge r-1$.
As $|Y'|=d_G(x')-1=r-1$,
 $|N_G(S)\cap Z|<|S|$
implies that
$S\ne Y'$ and thus $|S|\le r-2$.
For each $y\in Y'$,
$|N_G(y)\cap Z|\ge r-|\{u,x,x'\}|=r-3$,
implying that $|N_G(S)\cap Z|\ge r-3$.
Thus, $|N_G(S)\cap Z|<|S|$
implies that $|S|\ge r-2$.
Hence $|S|=r-2$, and $|N_G(S)\cap Z|=r-3$.

Let $Z':=N_G(S)\cap Z$. Then $|Z'|=r-3$.
For each $y\in S$,
$N_G(y)\subseteq \{u,x,x'\}\cup Z'$.
Thus,  for each $y\in S$,
since $d_G(y)=r$ and $|Z'|=r-3$,
we have $N_G(y)=\{u,x,x'\}\cup Z'$,
implying that  $y\in N_G(u)$
and thus
$S\subseteq N_{G}(u)\cap Y=Y_0$.
It follows that
$r-2=|S|\le |Y_0|=\frac r2$,
implying that $r\le 4$,
 contradicting
$r\ge 6$.
 Claim 4 holds.

\inclaim $G[X_0\cup Y_0]$ is not a complete bipartite graph.

Suppose that $G[X_0\cup Y_0]$ is a complete bipartite graph.
Then, for each $x_0\in X_0$,
$N_G(x_0)=\{u\}\cup Y_0\cup Y_1$ by Claim 4.
Since $|Y_0|=\frac r2$ by Claim 1,
we have  $|Y_1|=\frac r2-1$.

Since $N_G(u)=X_0\cup Y_0$,
$u$ is isolated by a perfect matching
$M_0$ of $G[X_0\cup Y_0]$.
Clearly, $M_0$ is a minimal
matching of $G$ which isolates
$u$.
By Corollary \ref{Krr},
$G_0:=G-(V(M_0)\cup \{u\})$
is isomorphic to
$K_{2t}$ or $K_{t,t}$
for some positive integer $t$.
Note that
$V(G_0)$ is a bipartite graph
with a bipartition
$(X_1\cup X_2, Y_1\cup Y_2)$,
as shown in
Figure~\ref{char-Fr-f2}.
Since $|X|\geq r$ and $|X_0|=\frac r2$, $|X_1\cup X_2|=|X|-|X_0|\ge \frac r2\ge 3$, implying that $G_0$ cannot be
a complete
graph.
It follows that $G_0$ is a complete bipartite graph with a bipartition
$(X_1\cup X_2, Y_1\cup Y_2)$,
as shown in
Figure~\ref{char-Fr-f2}.

Since
$G[X_1\cup X_2\cup Y_1\cup Y_2] $
(i.e., $G_0$) is a complete bipartite graph,
for any $x_1\in X_1$
we have
$Y_1\cup Y_2\subseteq N_G(x_1)$.
Then, by Claim 3,
$N_G(x_1)=Y_0\cup Y_1\cup Y_2$.
If there exists $x_2\in X_2$,
then,  we also have
$Y_1\cup Y_2\subseteq N_G(x_2)$.
By the assumption of $X_2$,
we have $N_G(x_2)\subseteq Y_1\cup Y_2$,
implying that
$Y_1\cup Y_2=N_G(x_2)$.
Hence
$|Y_1|+|Y_2|=d_G(x_2)=d_G(x_1)
=|Y_0|+|Y_1|+|Y_2|$,
implying that $|Y_0|=0$,
contradicting the fact that $|Y_0|=\frac r2\ge 3$.
Hence $X_2=\emptyset$.
Similarly, $Y_2=\emptyset$.

Now we know that $G[X_1\cup Y_1]$  is  a complete bipartite
graph.
By Claim 3,
$G[X_1\cup Y_0]$ is a complete bipartite graph.
Thus,
for each $x_1\in X_1$, we have $N_G(x_1)=Y_0\cup Y_1$.
Thus, $r=|Y_0|+|Y_1|$,
implying that $|Y_1|=r-|Y_0|=\frac r2$,
contradicting the conclusion that
$|Y_1|=\frac r2-1$
which is
obtained in the first paragraph of this claim's proof.
Claim 5 holds.

\inclaim $E(X_0, Y_0)\ne \emptyset$.

Suppose that $E(X_0, Y_0)= \emptyset$.
Then,
$N_G(x_0)=\{u\}\cup Y_1$
for any $x_0\in X_0$,
implying that $|Y_1|=r-1$.
Similarly, $|X_1|=r-1$.

Now we claim that $G[X_1\cup Y_1]$
is a complete bipartite graph.
Suppose that $x_1y_1\notin E(G)$,
where $x_1\in X_1$ and $y_1\in Y_1$.
Since $|X_1|=|Y_1|=r-1>\frac r2$ and
$G[X_0\cup Y_1]$ is a complete
bipartite graph by Claim 3.
Since $|Y_1|-|X_0|
=r-1-\frac r2\ge 2$,
$G[X_0\cup Y_1]$
has a matching
$M_1$ with $X_0\subseteq V(M_1)$ and $y_1\notin V(M_1)$.
Similarly,
$G[Y_0\cup X_1]$ has a matching
$M_2$ with
$Y_0\subseteq V(M_2)$
and $x_1\notin V(M_2)$.
Clearly, $M_1\cup M_2$ is a minimal matching of $G$ isolating $u$.
By Corollary \ref{Krr},
$G_0:=G-(V(M_1\cup M_2)\cup \{u\})$
is isomorphic to
$K_{2t}$ or $K_{t,t}$
for some positive integer $t$.
Note that the vertex set of $G_0$
can be partitioned into
two independent sets
$(X_1\setminus V(M_2))\cup X_2$
and $(Y_1\setminus V(M_1))\cup Y_2$.
Since $x_1\in X_1\setminus V(M_2)$
and $y_1\in Y_1\setminus V(M_1)$,
we have $x_1y_1\in E(G_0)$, contradicting the assumption that $x_1y_1\notin E(G)$.

Thus, $G[X_1\cup Y_1]$ is a complete bipartite graph.
By Claim 3,
$G[X_1\cup Y_0]$ is a complete bipartite graph.
Therefore, for any $x\in X_1$,
$Y_0\cup Y_1\subseteq N_G(x)$,
implying that
$d_{G}(x)\geq |Y_0|+|Y_1|=\frac{r}{2}+(r-1)>r$, a contradiction.
Claim 6 holds.

\inclaim $|X_1|=\frac r2$ and $|Y_1|=\frac r2$.

By Claim 5, $G[X_0\cup Y_0]$ is not a complete bipartite graph.
Assume that $x_0y_0\notin E(G)$,
where $x_0\in X_0$ and $y_0\in Y_0$.
Since $N_G(x_0)\subseteq \{u\}\cup (Y_0\setminus \{y_0\})\cup Y_1$ and $d_G(x_0)=r$,
we have
$r\le 1+(\frac r2-1)+|Y_1|$,
implying that $|Y_1|\ge \frac r2$.
Similarly, $|X_1|\ge \frac r2$.

By Claim 6, $E_G(X_0,Y_0)\ne \emptyset$.
Assume that $x_1y_1\in E(G)$, where
$x_1\in X_0$ and
$y_1\in Y_0$.
Let $Y'=N_G(x_1)\cap Y$.  Obviously, $y_1\in Y'$.
By Claim 3, $Y_1\subseteq Y'$.
By Claim 3 again,
$X_0\cup X_1\subseteq N_G(Y_1)\cup N_G(y_1)\subseteq N_G(Y')$.
By Claim 4, $|N_G(Y')\cap X|\le r$,
implying that $|X_0|+|X_1|
\le |N_G(Y')\cap X|\le r$.
Since $|X_0|=\frac r2$, we have
$|X_1|\le \frac r2$.
Due to the conclusion in the previous paragraph, we have $|X_1|=\frac r2$.
Similarly, $|Y_1|=\frac r2$.
The claim holds.

\inclaim $X_2=Y_2=\emptyset$.

Suppose that $X_2\ne \emptyset$.
By the definition of $X_2$,
for each $x\in X_2$,
we have
$N_G(x)\subseteq Y_1\cup Y_2$.
Since $|Y_1|=\frac r2$ by Claim 7,
$r=d_G(x)\le |Y_1|+|Y_2|=\frac r2+|Y_2|$ and thus $|Y_2|\ge \frac r2$.
Obviously, $|X_2|=|Y_2|\ge \frac r2$.

By Claim 5,
there exist $x_0\in X_0$ and $y_0\in Y_0$ such that $x_0y_0\notin E(G)$.
Since $G$ is connected,
we have $E_G(X_2\cup Y_2, X_1\cup Y_1)\ne \emptyset$.
Without loss of generality,
assume that $x'y'\in E(G)$,
where $x'\in X_2$ and $y'\in Y_1$.

Since $G[X_0\cup Y_1]$ is
a complete
bipartite graph by Claim 3,
$G[(X_0\cup Y_1)]- \{x_0,y'\}$ has a perfect matching $M'$.
Similarly,
$G[Y_0\cup X_1]$ has a matching
$M''$ with
$V(M'')\cap Y_0=Y_0\setminus \{y_0, y_1\}$, where $y_1\in Y_0\setminus \{y_0\}$.
Thus, $M_0:=\{uy_1,x'y'\}\cup M'\cup M''$
is a matching of $G$ isolating $x_0$.
Since $N_{G}(x_0)\subseteq \{u\}\cup Y_0\cup Y_1$ and $|Y_0|= |Y_1|=\frac{r}{2}$, then $y_0$ is the only vertex in $Y_0$ that is not adjacent to $x_0$. It follows that $M_0$ is a minimal matching of $G$ isolating $x_0$.
Observe that
the vertex set of
$G_0:=G-(V(M_0)\cup \{x_0\})$
can be partitioned into independent sets
$(X_2\setminus \{x'\})\cup (X_1\setminus V(M''))$
and $\{y_0\}\cup Y_2$.
By Corollary \ref{Krr},
$G_0\cong K_{t,t}$
where
$t=|Y_2|+1\ge \frac r2+1\ge 4$,
contradicting the fact that
$y_0$ is a vertex in $G_0$ with
$d_{G_0}(y_0)=|X_1\setminus V(M'')|=2$.
The claim holds.

By Claim 8,
$X=X_0\cup X_1$ and
$Y=Y_0\cup Y_1$.
By Claims 1 and 7,
$|X_0|=|X_1|=|Y_0|=|Y_1|=\frac r2$.
Then, for each $x\in X_1$, we have
$r=d_G(x)\le |Y_0|+|Y_1|=r$,
implying that $Y_0\cup Y_1
\subseteq N_G(x)$.
Similarly, $X_0\cup X_1\subseteq N_G(y)$ for each $y\in Y_1$.
By Claim 3, we have $|N_G(x)\cap Y_0|=\frac r2-1$ for each $x\in X_0$
and $|N_G(y)\cap X_0|=\frac r2-1$ for each $y\in Y_0$.
Hence $G\cong F_r$.
\proofend

We are now ready to prove the main result.

\vspace{0.3 cm}

\setcounter{countclaim}{0}

\noindent {\it Proof of Theorem~\ref{main1}}:
Let $G=(V,E)\in \G$.
By Proposition \ref{M22}
and Lemma~\ref{basic-2} (iv),
$M^v_{2,2}$ contains exactly one
edge, say $uu'$. Assume that
$|N_G(u)\cap W|\ge 2$.
Let $I':=W\cup \{u'\}$.
We will complete the proof by
establishing two claims below.

\inclaim
$G-(I\cup \{u\})$ does not have
any odd clique $X$
with the properties that $|X|\ge 3$,
$|X\cap I'|=1$ and
$N_G(X)\subseteq \{u\}\cup X\cup I$.

By Lemma~\ref{yT''} (i),
$I'$ is an independent set of $G$.
By Lemma~\ref{basic-2} (iv),
$u$ is the only vertex in
$V\setminus (I\cup W)$
with the property that
$|N_G(u)\cap W|\ge 2$.
Hence, $|N_{G}(h)\cap I'|\leq 1$ for each $h\in V\setminus(I\cup I'\cup\{u\})$.
Also note that $\matt\setminus \{uu'\}$
is a perfect matching of the subgraph
$G-(I\cup I'\cup \{u\})$.
By Proposition~\ref{move-T},
Claim 1 holds.

\inclaim $T''=W'=\emptyset$.

If $T''\ne \emptyset$, then,
Lemma~\ref{yT''} (ii) implies that
$X:=T''\cup \{u'\}$ is an odd clique
with the property that $|X|\ge 3$,
$X\cap I'=\{u'\}$ and
$N_G(X)\subseteq I\cup X\cup \{u\}$,
contradicting Claim 1.
Thus, $T''=\emptyset$.

Recall that $W'$ denotes the set of vertices $w\in W$ for which $|T_w|\ge 3$.
Suppose that $W'\ne \emptyset$ and $w'\in W'$.
By Lemma~\ref{yT''} (iii),
$X:=T_{w'}$ is an odd clique.
As $w'\in W'$, we have
$|X|\ge 3$.
By Lemma~\ref{yT''} (iii),
$N_G(u')\cap W'=\emptyset$,
and thus $X\cap I'=\{w'\}$.

By Lemma~\ref{basic-2} (ii),
$N_G(X)\cap W\subseteq \{w'\}
\subseteq X$.
If $N_G(X)\not \subseteq I\cup X\cup \{u\}$,
then $T''=\emptyset$ implies that
$N_G(X)\cap T_{w}\ne \emptyset$ for some $w\in W'\setminus \{w'\}$.
By Lemma~\ref{yT''} (iii),
$T_{w}$ is an odd clique.
By Lemma~\ref{basic-2} (ii) again,
$w\notin N_G(X)$
and $w'\notin N_G(T_w)$.
Then,
$N_G(X)\cap T_{w}\ne \emptyset$
implies the existence of an
edge $vv'$ in $G$, where $v\in T_{w}\setminus \{w\}$
and $v'\in X\setminus \{w'\}$.
Since both $X$ and $T_w$ are
cliques, $G[X\cup T_w]$ contains
an $M^v_2$-augmenting path
between $w$ and $w'$ which
contains edge $vv'$,
contradicting
Lemma~\ref{basic-2} (i).
Hence $N_G(X)\subseteq I\cup X\cup \{u\}$.

Now it has been shown that $X$ is an odd clique of $G-(I\cup \{u\})$ with the properties that $|X|\ge 3$, $X\cap I'=\{w'\}$ and $N_G(X)\subseteq I\cup X\cup \{u\}$.
However, by
Claim 1, such a set $X$ does not exist.
Hence $W'=\emptyset$ and
Claim 2 holds.

Finally, we are going to complete the proof.  By Claim 2,
$T''=\emptyset$ and $W'=\emptyset$.
Thus, $M^v_{2,j}=\emptyset$
for both $j=0,1$.
Then, $M^v_{2,2}=\{uu'\}$ implies that  $V=I\cup W\cup \{u,u'\}=I\cup I'\cup \{u\}$.
By Lemma~\ref{yT''} (i),
$I'$ is an independent set.
Thus,
$V\setminus \{u\}$ can
be partitioned into two
independent sets
$I$ and $I'$
with $|I|=|I'|$.
By Lemma~\ref{char-Fr}, we have
$G\cong F_r$.
Hence Theorem~\ref{main1}  holds.
\proofend

\section{Remarks and Problems}

Theorem~\ref{main1}
contributes the final piece of solutions
to the characterization
of connected equimatchable $r$-regular graphs
(see a schematic overview
in Figure~\ref{lastpiece}).
On the study of
equimatchable graphs which may not be regular,
B\"{u}y\"{u}kolak et al. \cite{YB}
 provided a complete structural characterization of equimatchable graphs without triangles, by identifying the equimatchable triangle-free graph families.

However, the structure of equimatchable graphs with triangles has not been fully characterized yet. $F_r$ is a special class of equimatchable graphs with triangles, which is constructed from complete bipartite graphs.
It is worth exploring whether equimatchable graphs with triangles can be constructed from triangle-free equimatchable graphs.

\section*{Declarations}

\noindent 
{\bf Funding Declaration}: 
This research is supported by
Natural Science Foundation of China  
(No.12101203, 12371340)
and Natural Science Foundation of Henan Province (No.252300421489).

\noindent {\bf 
Consent to Participate declaration}: not applicable

\noindent {\bf 
	Consent to Publish declaration}: not applicable

\noindent {\bf 
Ethics declaration}: not applicable

\end{document}